\documentclass[
 hidempi, 
]{mpi2015-cscpreprint}

\usepackage{graphicx} 
\usepackage{longtable}
\usepackage{lscape}
\usepackage{a4wide}
\usepackage{tabularray}
\usepackage{xcolor}
\usepackage{here}
\usepackage{amsmath}
\usepackage{amssymb}
\usepackage{natbib}
\usepackage{pgffor}
\usepackage{color}
\usepackage{cleveref}
\usepackage{authblk}
\usepackage{subcaption}
\usepackage{multirow}
\usepackage{arydshln}
\usepackage{graphbox}
\usepackage{colortbl}
\usepackage{amsfonts}
\usepackage{epsfig}
\usepackage{colortbl}
\usepackage{comment}
\usepackage{mathdots}
\usepackage{mathrsfs}
\usepackage{epstopdf}
\usepackage{makeidx} 
\usepackage{ifpdf}
\usepackage{arydshln}
\usepackage{tikz,tikz-3dplot}
\usetikzlibrary{arrows,decorations.pathmorphing,backgrounds,fit,positioning,shapes.symbols,chains,mindmap,trees,calc}
\usepackage{xfp}

\definecolor{codegreen}{rgb}{0,0.6,0}
\definecolor{codegray}{rgb}{0.5,0.5,0.5}
\definecolor{codepurple}{rgb}{0.58,0,0.82}
\definecolor{backcolour}{rgb}{0.95,0.95,0.92}

\title{Rational approximation for Zolotarev sign and ratio problems}

\author[$\S$]{C. Poussot-Vassal}
\affil[$\S$]{DTIS, ONERA, Universit\'e de Toulouse, 31000, Toulouse, France.\authorcr
	\email{charles.poussot-vassal@onera.fr}, \orcid{0000-0001-9106-1893}}

\author[$\dag$]{I. V. Gosea}
\affil[$\dag$]{Max Planck Institute for Dynamics of Complex Technical Systems, CSC Group \\

	Sandtorstr. 1, 39106 Magdeburg, Germany.\authorcr
	\email{gosea@mpi-magdeburg.mpg.de}, \orcid{0000-0003-3580-4116}}

\author[$\ast$]{A. C. Antoulas}
\affil[$\ast$]{Department of Electrical and Computer Engineering, Rice University, Houston, TX, USA.\authorcr
	\email{aca@rice.edu}}

\abstract{We study algorithms for Zolotarev (3rd and 4th) problems rational approximation. First, we show that the Loewner framework (LF) is appropriate to rapidly and with no iterations, approximate the Zolotarev problems by compressing the (numerous) interpolation conditions. Second, we compare the approximation properties (e.g. coefficients and poles) of LF with the standard AAA, AAA-sign and its AAA-Lawson variants. We concentrate the study on a canonical example, namely the symmetric two-circles one, for which the optimal solution is well documented. For this case, we highlight the numerical robustness of LF and its ability to recover the structure of the optimal solution. Additional non-trivial geometries are also reported, emphasizing that LF is fast, reliable, and yields accurate approximants with no user intervention. 
}


\novelty{We show that the Loewner Framework (LF) (i) is a viable candidate for approximating the third and fourth Zolotarev problems; (ii) is accurate, simpler, faster and numerically more robust than AAA; and (iii) in some specific cases, recovers the mathematical structure of the optimal solution.
}

\definecolor{blue}{RGB}{16,97,169} 
\definecolor{grey}{RGB}{88, 102, 110} 
\definecolor{green}{RGB}{65,209,204} 
\definecolor{orange}{RGB}{224, 131, 0} 
\definecolor{pink}{RGB}{255,105,180} 

\def\IR{{\mathbb R}}
\def\IC{{\mathbb C}}

\def\IS{{\mathbb S}}
\def\IL{{\mathbb L}}
\def\sIL{{\IS}}

\newcommand{\bA}{{\mathbf A}}
\newcommand{\bB}{{\mathbf B}}
\newcommand{\bC}{{\mathbf C}}

\newcommand{\bE}{{\mathbf E}}
\newcommand{\bF}{{\mathbf F}}

\newcommand{\bS}{{\mathbf S}}
\newcommand{\bY}{{\mathbf Y}}

\newcommand{\bW}{{\mathbf W}}
\newcommand{\bR}{{\mathbf R}}

\newcommand{\bX}{{\mathbf X}}

\newcommand{\bV}{{\mathbf V}}

\newcommand{\bh}{{\mathbf h}}

\newcommand{\cT}{ {\cal T} }

\newcommand{\varIC}{z}
\newcommand{\rloe}{\bR_{{\mathrm{Loew}}}}
\newcommand{\raaa}{\bR_{{\mathrm{AAA}}}}
\newcommand{\chebfun}{{\texttt{chebfun}}}
\def\signEF{\mathrm{sign}_{E/F}}

\newtheorem{theorem}{Theorem}
\newtheorem{lemma}{Lemma}
\newtheorem{remark}{Remark}

\newtheorem{problem}{Problem}

\definecolor{blue}{RGB}{16,97,169} 
\definecolor{grey}{RGB}{88, 102, 110} 
\definecolor{green}{RGB}{65,209,204} 
\definecolor{orange}{RGB}{224, 131, 0} 
\definecolor{pink}{RGB}{255,105,180} 
\definecolor{colhl}{RGB}{152,251,152}
\definecolor{colhl}{RGB}{211,211,211}

\def\CAS{10}
\def\NARCH{6}

\begin{document}

\maketitle

\addtocontents{toc}{\setcounter{tocdepth}{2}} 
\tableofcontents

\begin{center}
\vspace{1cm}
\includegraphics[width=.45\textwidth]{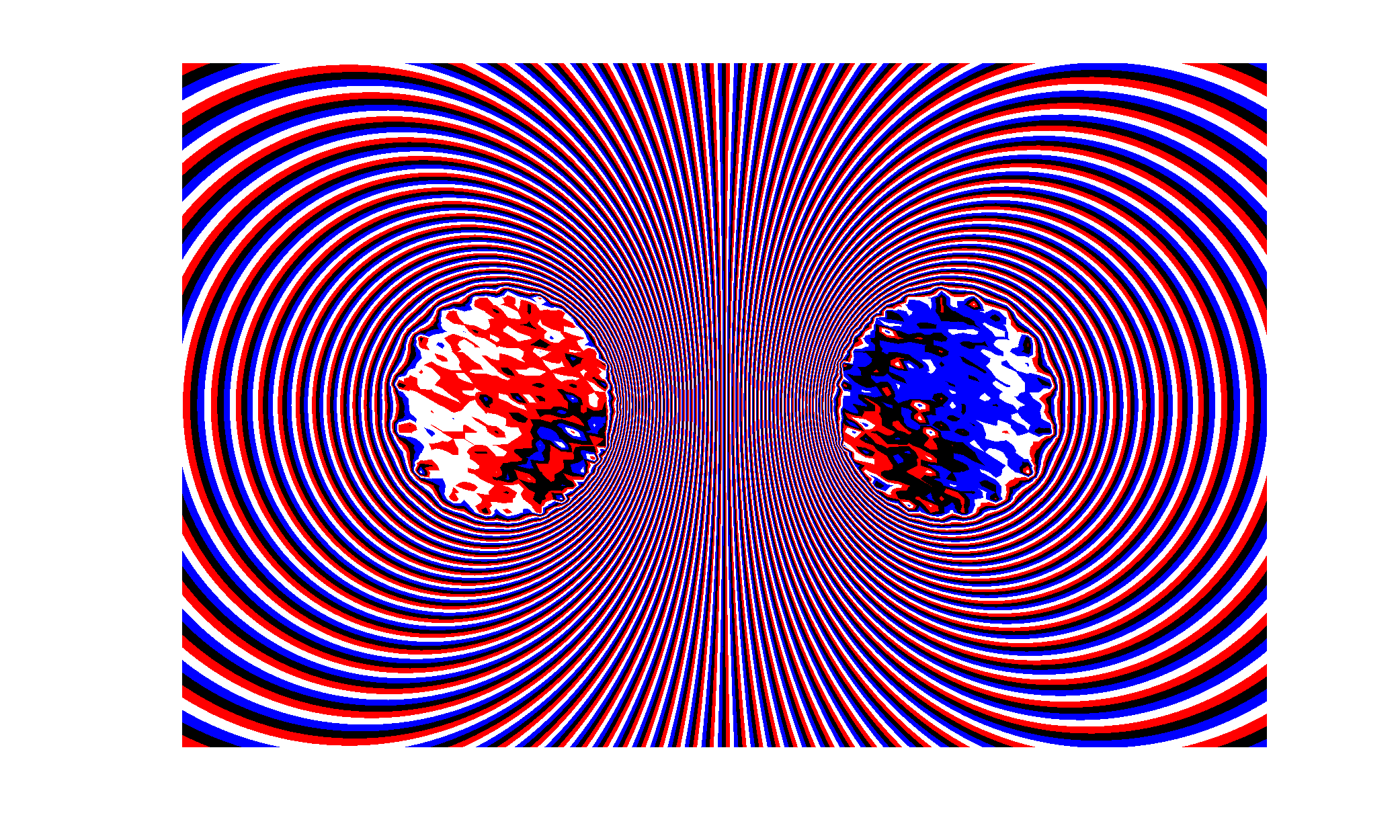}
\includegraphics[width=.45\textwidth]{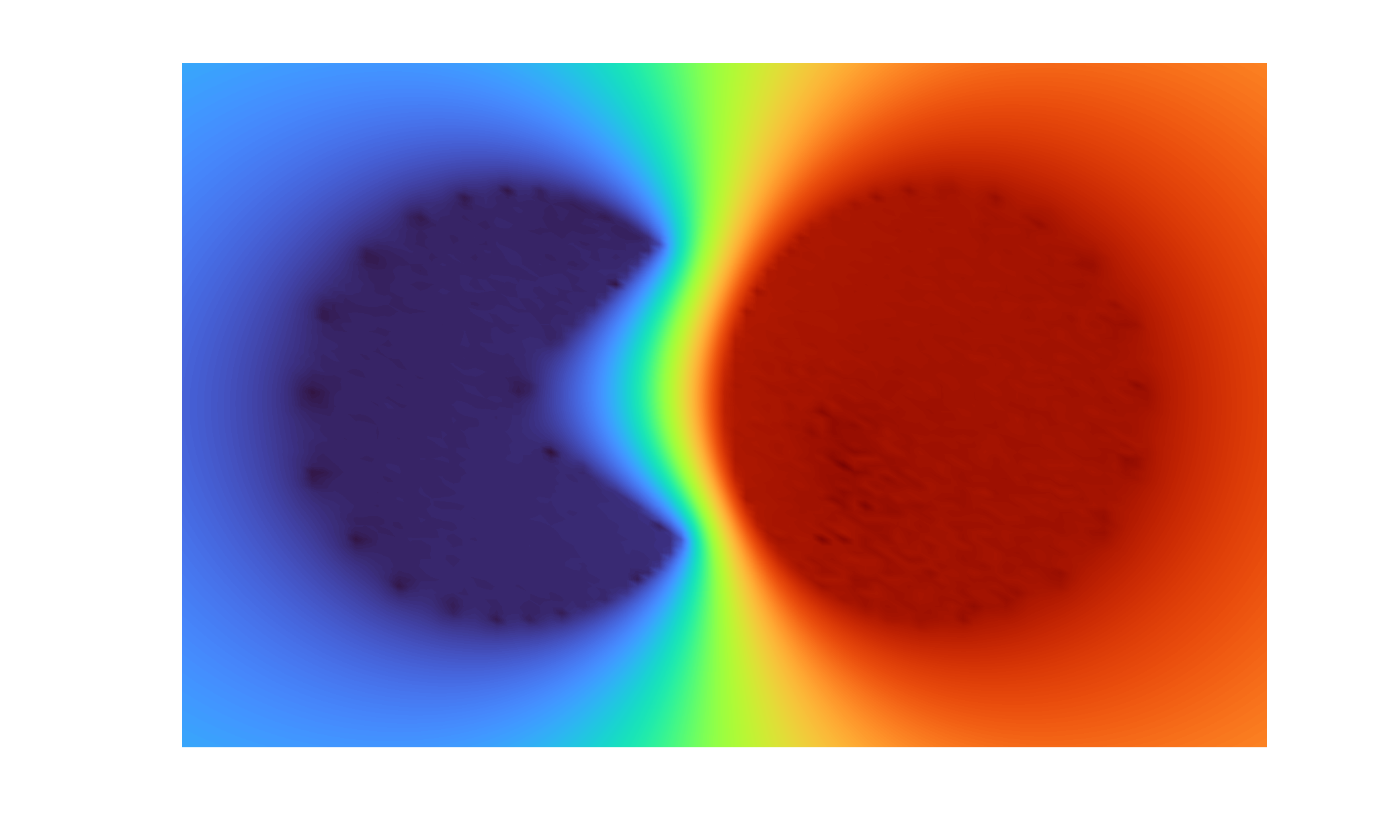}\\
Artistic views of the third Zoloratev problem applied to \\ the two symmetric circles (left) and Pac Man (left) topologies.
\end{center}

\section{Introduction and contributions}
\label{sec:intro}

\paragraph{About rational approximation.}

Approximation theory is a classic,  well-established field of mathematics, split into a wide variety of sub-fields, out of which we are particularly interested in approximation by rational functions. This has a multitude of applications in areas of applied sciences. Some examples include computational fluid dynamics, solution and control of partial differential equations, data compression, electronic structure computation, systems and control theory (model order reduction of dynamical systems).  For more details of the many aspects of approximation theory, we refer the reader to the classical books in \cite{Da75,Po81}, and to the more recent books \cite{Tr13,Is18} that also put emphasis on practical issues such as software implementation of the methods and applications, i.e., on data analysis for the latter.

Since rational functions are ratios of two polynomials, they have proven to be more versatile in approximating complex functions than simple polynomials. Furthermore, rational functions are extensively used in signal processing, systems, and control theory since the input-output behavior of linear dynamical systems is described in the frequency domain by the so-called transfer function (which is a rational function). Finally, the location of the poles of transfer functions determines important system properties such as asymptotic stability, transient behavior, damping, and harmonics.

\paragraph{Preliminaries on the Zoloratev 3rd and 4th problems and the data-driven framework.}

We are concerned with rational approximation, i.e.,  approximation of any type of functions, by means of rational ones. Here, we consider sign and ratio Zolotarev functions whose solution is intertwined. More precisely, the original function to be approximated falls into the category of a special sign function, defined on the complex domain, i.e., one that is defined on a reunion of two domains $E$ and $F$, for which the value is constant on each, e.g., equal to $\pm 1$ (hence the name Zolotarev sign problem, Z4, \cref{pb:z4}). To find the solution to the Zolotarev ratio function, one may appropriately transform the previously obtained solution. What is obtained is a function that is large on one domain ($F$), and small on the other ($E$), and thus, their ratio is minimized (hence the name Zolotarev ratio problem, Z3, \cref{pb:z3}). 

Moreover, the methods we use to solve Zolotarev sign problems do not require access to the exact closed-form of the original function, but only to evaluations of it on a particular domain. Hence, these methods are \emph{data-driven}.

In this work, we apply two established data-driven rational approximation methods: the Loewner Framework (LF) introduced in \cite{ajm07}, and the Adaptive Antoulas-Anderson (AAA) algorithm introduced in \cite{nst18}, and further developed in \cite{nakatsukasa2020algorithm,trefethen2024computation}. 
The methods mentioned above have some things in common: they both compute approximants in the form of rational functions, they exclusively use measurements of the original function (and not the function itself), and they involve, in one way or another, Loewner matrices.

\subsection{Contributions}
We summarize below the main contributions of this paper, which were confirmed through a series of comprehensive numerical experiments:
\begin{itemize}
    \item[(i)] we show that the LF solves Zolotarev problems by compressing the (many) number of  interpolation points. This approach does not need neither any iterations nor user intervention, and yields solutions quite close to the optimal ones in a very low computational time. Consequently, the LF constitutes a valid alternative for addressing these problems;
    \item[(ii)] we show that, in some specific cases (notably the symmetric two circles one), LF recovers the structure of the optimal solution as well as its symmetry property (e.g., the distribution of eigenvalues and zeros, the alternance of polynomials coefficients, realness, etc.), whereas the other methods tend to add spurious and oddly distributed poles and zeros, and non-trivial artifacts that can't be easily explained. In addition, thanks to an extensive numerical experimentation, we demonstrate the numerical robustness in the results obtained with LF compared to AAA\footnote{Up to the author's knowledge, this property has not been highlighted in previous works.};
    \item[(iii)] we conduct an extensive numerical study\footnote{Results can be re-computed using the code available at \url{https://github.com/cpoussot/zolotarev}.} and comparison of the performance of different methods\footnote{The AAA algorithm uses the code available at \url{https://github.com/chebfun/chebfun}.} w.r.t. to computing time, accuracy of fit and interpretability, for approximating several sign functions defined on various domains showing, among others, the computational advantage of LF.
\end{itemize}

\subsection{Paper structure}
This paper is organized as follows: after the introduction in \cref{sec:intro}, in \cref{sec:zolo}, the Zolotarev sign and ratio problems are presented. Then, \cref{sec:ratapp} describes the main rational approximation tools to be used later on, e.g., the Loewner Framework in \cref{sec:Loew} and the AAA algorithm in \cref{sec:AAA}. We keep the exposition brief, since both are established methods, and we refer to the original manuscripts for details. Then, in \cref{sec:1a}, the special symmetric two circles topology is treated (i.e., for which the $E$ and $F$ domains are two symmetric circles). For this configuration, a detailed analysis of the results obtained by LF and AAA is discussed. The latter is not restricted to the approximation metric but also includes the values of the obtained numerators and denominators polynomials, poles and zeros, etc. (more details are available in \cref{app:1a}). Then,  in \cref{sec:num}, we conduct an exhaustive numerical experiments for various test cases, including the ones in \cite{trefethen2024computation} (more details are available in \cref{app:all}). Finally, applications of Zolotarev problems are presented in  \cref{app:applications}. Conclusions are discussed in \cref{sec:conclusion}.

\section{The Zolotarev problems}\label{sec:zolo}

In \cite{trefethen2024computation}, a detailed historical account is presented, which is summarized next. Additionally, we refer the reader to the review paper \cite{nakatsukasa2016computing} for further considerations.

\subsection{Reminder and general procedure}

For integers $m,n\ge 0$, let $R_{m,n}$ (or simply $R_n$ if $m=n$) be the set of rational functions with numerator of degree $m$ and denominator of degree $n$, which means that any rational function $r\in R_{m,n}$ can be written as $p(\varIC)/q(\varIC)$ for some polynomials $p(\varIC)$ and $q(\varIC)$ of degrees at most $m$ and $n$, respectively. The problem  studied here is to find $r_{m,n}^*\in R_{m,n}$ such that
\begin{equation}
r_{m,n}^* := \arg \min_{r\in R_{m,n}} {\max_{\varIC\in E} |r(\varIC)| \over \min_{\varIC\in F} |r(\varIC)|}.
\label{eq:unscaled}
\end{equation}
Through normalization, i.e., assuming that $\min_{\varIC\in F} |r(\varIC)| = 1$, we end up with the problem of minimizing $|r(\varIC)|$ solely over set $E$. This is known in the literature as the third Zolotarev problem, or sometimes, also the {\em Zolotarev ratio problem}.

\begin{problem}[Zolotarev ratio problem - Z3]\label{pb:z3}
Find $r_{m,n}^*\in R_{m,n}$, satisfying condition $\min_{\varIC\in F} |r_{m,n}^*(\varIC)| = 1$, and that attains the minimum (where $\|\cdot\|_E^{}$ denotes the supremum norm over $E$)
\begin{equation}
\sigma_n = \min_{r \in R_{m,n}} \|r\|_E^{}.
\label{eq:Z3}
\end{equation}
\end{problem}

As mentioned in \cite{trefethen2024computation}, it is known that a solution $r_{m,n}^*$ to \cref{eq:Z3} exists, though it need not be unique, and that $\sigma_n$ satisfies $0 < \sigma_n \le 1$.

Then, the {\em sign function} relative to $E$ and $F$ is of relevance in what follows, and is defined as follows:
\begin{equation}
\signEF(\varIC) =
\begin{cases}
-1 & \varIC\in E, \\ +1 & \varIC\in F.
\end{cases}
\end{equation}
It is to be noted that for $z\in\IC\backslash \{E\cup F\}$, $\signEF(\varIC)$ remains undefined.

Next, we introduce the problem of rational minimax approximation of $\signEF$ over $E$ and $F$, or problem Z4: 

\begin{problem}[Zolotarev sign problem - Z4]\label{pb:z4}
Find $\hat r_{m,n}\in R_{m,n}$ that attains the minimum (where $\|\cdot\|_{E\cup F}$ denotes the supremum norm over $E\cup F$)
\begin{equation}
\tau_n = \min_{r\in R_{m,n}} \|r-\signEF\|_{E\cup F}^{}.
\label{eq:Z4}
\end{equation}
\end{problem}

As mentioned in \cite{trefethen2024computation}, the proof showing that the problems Z3 and Z4 are equivalent in the case for which $E$ and $F$ are real intervals is due to the work of \cite{akhiezer1990elements}. Then, in \cite{istace1995third}, it was shown that the equivalence generalizes to the complex case.  We reproduce (with appropriate notation changes) in what follows their theorem.

\begin{theorem}
\label{thm:thm1}
Every solution $r_{m,n}^*$ of Problem Z3 is related to a solution 
$\hat r_{m,n}$ of Problem Z4 by
\vspace{-2mm}
\begin{equation}
\hat r_{m,n} (\varIC) = {1-\sigma_n\over 1+\sigma_n} \,{r_{m,n}^*(\varIC) - \sqrt{\sigma_n}
\over r_{m,n}^*(\varIC) + \sqrt{\sigma_n}}, \quad
r_{m,n}^*(\varIC) = \sqrt{\sigma_n} \,
{(1-\sigma_n)/(1+\sigma_n) + \hat r_{m,n}(\varIC) \over
(1-\sigma_n)/(1+\sigma_n) - \hat r_{m,n}(\varIC)}.
\label{eq:relation}
\end{equation}
The minimal values of the two problems satisfy
\vspace{-2mm}
\begin{equation}
\tau_n = {2\sqrt\sigma_n\over 1 + \sigma_n}, \quad
\sigma_n = \left({\tau_n\over 1 + \sqrt{1-\tau_n^2}\kern .7pt} \right)^2,
\label{eq:minvals}
\end{equation}
and the set of extremal points $M = \{z\in E\cup F, ~|\hat r_{m,n}(\varIC) - \signEF(\varIC)| = \tau_n\}$ is the union of  $M_1 = \{z\in E, ~|r_{m,n}^*(\varIC)| = \sigma_n\}$ and $M_2 = \{z\in F, ~|r_{m,n}^*(\varIC)| = 1\}$.
\end{theorem}

\subsection{Approach paved in  \cite{trefethen2024computation} and discussions}

The proposed algorithms for solving \cref{pb:z3} (Z3) actually consist of solving \cref{pb:z4} (Z4) first, and then transforming from $\hat r_{m,n}$ and $\tau_n$ to $r_{m,n}^*$ and $\sigma_n$. Following \cref{thm:thm1}, the authors of \cite{trefethen2024computation} approach Problem Z3 in the following two-step procedure: (i) solve Problem Z4 using a rational approximation method; then (ii) convert to a solution of Problem Z3 using \cref{eq:relation} and \cref{eq:minvals}.

The main challenge of the algorithm proposed in \cite{trefethen2024computation} is given by step (i), i.e., the computation of the degrees $(m,n)$ rational best approximation $\hat r_{m,n}$ to the sign function $\signEF$ defined by the sets $E$ and $F$.  In \cite{trefethen2024computation}, the proposed procedure is to find $\hat r_{m,n}$ by means of the AAA approximation technique.  However, this has proven quite challenging. The first difficulty reported in \cite{trefethen2024computation} is that AAA encounters particular challenges when applied to sign functions. The second is the importance of having not just good rational approximations but nearly optimal ones, requiring the use of the AAA-sign and AAA-Lawson algorithm \cite{nakatsukasa2020algorithm}, which had a known problem of non-convergence in certain cases. These difficulties and their proposed solutions are discussed in \cite[\S 3.1 and \S 3.2]{trefethen2024computation}. Basically, the Lawson iteration "drives" the standard AAA approximant close to the optimal approximation value by performing several least squares steps involving  the weights. Although seemingly effective, this approach has proven to be quite costly, mostly in terms of running time, but also in terms of the structure of the resulting polynomials, as will become evident in the next sections. This is particularly relevant when choosing a high degree $(m,n)$ for the approximant. In addition, as highlighted in the coming sections, AAA (and its extensions) results vary according to the computer platform used.

\section{The main approximation methods}\label{sec:ratapp}

\subsection{The Loewner framework}\label{sec:Loew}

The Loewner Framework (LF) is a data-driven model identification, approximation, and reduction tool introduced in \cite{ajm07}. It provides a simple solution through rational approximation by means of interpolation. For more details, we refer the reader to the extensive tutorial/handbook papers in \cite{morKarGA19a,morGosPA22,ali17}, and to the papers \cite{morIonA14,AGPV:2025,morAntGI16,morKarGGetal25} that address extensions of the LF to some special classes of parametric, nonlinear or hybrid dynamical systems, and to  the recent book of interpolation and model reduction \cite[chapter 4]{abg_book}. 



Given $N>0$ distinct sampling points $\cT = \{\tau_1, \tau_2, \ldots,\tau_N\}$ together with $N$ function evaluations at these points denoted with $\{f_1, f_2, \ldots,f_N\}$, let $\mathfrak{D} = \{(\tau_\ell; \  f_\ell) \vert \ell =1,\cdots,N\}$ to be the problem data. Assume that $N=k+q$, where $k$ and $q$ are positive integers. The given data $\mathfrak{D}$ is first partitioned into the following two disjoint subsets,
\begin{equation}\label{eq:data_Loew}
\mathfrak{D} = \mathfrak{L} \cup \mathfrak{R}, \ \text{where} \ \begin{cases} {\text{right data}}: \mathfrak{R} = \{(\lambda_i; \  w_i) \vert i=1,\cdots,q
\}, \\
 {\text{left data}}: \ \ \mathfrak{L} =
\{(\mu_j; \ v_j) \vert j=1,\cdots,k\}. \end{cases}
\end{equation} 
The elements of the subsets $\mathfrak{L}$ and $\mathfrak{R}$ are as follows (for $1 \leq i \leq q$ and $1 \leq j \leq k$):
\begin{enumerate}
    \item  $\lambda_i$, $\mu_j$ $\in\IC$ are the right, and, respectively, left \emph{sample/interpolation points}.
	\item  $w_i$, $v_j$ $\in\IC$ are the right, and, respectively, left \emph{sample values}.
\end{enumerate}
The problem can be formulated as follows: find a rational function $\rloe(\varIC)$ of order $r \ll k$, such that the following interpolation conditions are (approximately) fulfilled:
\begin{equation} \label{eq:interp_cond}
\rloe(\lambda_i)=w_i,~~~\rloe(\mu_j)=v_j.
\end{equation}


The next step is to arrange data \cref{eq:data_Loew} into matrix format. The \textbf{Loewner matrix} $\IL \in\IC^{k\times q}$ and the \textbf{shifted Loewner matrix} $\sIL \in\IC^{k\times q}$ are defined as
\begin{equation} \label{eq:Loew_mat}
\IL=\left[\begin{array}{ccc}
\frac{v_1-w_1}{\mu_1-\lambda_1} & \cdots &
\frac{v_1-w_q}{\mu_1-\lambda_q} \\
\vdots & \ddots & \vdots \\
\frac{v_k-w_1}{\mu_k-\lambda_1} & \cdots &
\frac{v_k-w_q}{\mu_k-\lambda_q} \\
\end{array}\right], \ \ 
\sIL=\left[\begin{array}{ccc}
\frac{\mu_1v_1-w_1\lambda_1}{\mu_1-\lambda_1} & \cdots &
\frac{\mu_1v_1-w_q \lambda_q}{\mu_1-\lambda_q} \\
\vdots & \ddots & \vdots \\
\frac{\mu_kv_k-w_1 \lambda_1}{\mu_k-\lambda_1} & \cdots &
\frac{\mu_kv_k-w_q\lambda_q}{\mu_k-\lambda_q} \\
\end{array}\right],
\end{equation}
while the data vectors $\bV\in \IR^k, \bW^\top \in \IR^q$ are introduced as
\begin{equation} \label{eq:VW_vec}
\bV=\left[\!\begin{array}{cccc} v_1 \  \ v_2 \  \ \ldots \  \ v_k \end{array}\!\right]^\top \in \IC^k, \ \ \ \bW=[w_1~~w_2~~\cdots~~w_q] \in \IC^q.
\end{equation}

\begin{lemma}
	If $k=q$ the matrix pencil $(\sIL,\,\IL)$ is regular, then one can directly compute the  interpolation function as below:
	\begin{equation}
	\rloe(\varIC)=\bW(\sIL-\varIC \IL)^{-1}\bV.
	\end{equation}
\end{lemma}

In applications (see, e.g., \cite{morGosPA22}), the pencil $(\sIL,\,\IL)$ is often singular. In these cases, one needs to perform an SVD of the Loewner pencil, constructed from the row- and column-wise concatenation of matrices $\IL$ and $\sIL$. By choosing a truncation value $r$, we can write the truncated SVDs as follows (with $\bX_r = \bX_r^{(1)} \in \IC^{k \times r}$, $\bY_r = \bY_r^{(2)} \in \IC^{q \times r}$)
\begin{equation}\label{eq:loew-rank}
\begin{bmatrix} \IL & \sIL \end{bmatrix} \approx \bX_r^{(1)} \bS_r^{(1)} {\bY_r^{(1)}}^* ,~
\begin{bmatrix} \IL \\ \sIL \end{bmatrix}  \approx \bX_r^{(2)} \bS_r^{(2)} {\bY_r^{(2)}}^*. 
\end{equation}	
\begin{lemma}			
	The rational function corresponding to the reduced-order Loewner model is derived as:
    \begin{align*}
	\rloe(\varIC) = \hat{\bC}_{\mathrm{Loew}} (\varIC \hat{\bE}_{\mathrm{Loew}}-\hat{\bA}_{\mathrm{Loew}})^{-1} \hat{\bB}_{\mathrm{Loew}} 
    =  \bW \bY_r \Big{(} \bX_r^* (\sIL- \varIC \IL) \bY_r \Big{)}^{-1} \bX_r^* \bV,
	\end{align*}
    where the system matrices corresponding to a reduced Loewner model are:
    	\begin{equation} \label{eq:Loew_red_mat}
	\hat{\bE}_{\mathrm{Loew}} = -\bX_r^*\IL \bY_r, \ \  \hat{\bA}_{\mathrm{Loew}} = -\bX_r^*\sIL \bY_r, \ \
	\hat{\bB}_{\mathrm{Loew}} = \bX_r^*\bV, \ \  \hat{\bC}_{\mathrm{Loew}} = \bW \bY_r.
	\end{equation}
\end{lemma}

\subsection{The AAA algorithm (and variants)}\label{sec:AAA}

The AAA (Adaptive Antoulas-Anderson) algorithm proposed in \cite{nst18} represents an adaptive extension of the method originally introduced by A.C. Antoulas and B.D.O. Anderson in \cite{aa86}. The AAA algorithm is a rational approximation data-based method that makes use of the barycentric form of rational functions and uses an iteration to add more interpolation points (in a greedy fashion) until a certain level of approximation is attained. As exposed in \cref{sec:Loew}, it only requires as input evaluations of the function on a set of points denoted with $\cT$. We refer to \cite{AAAreview} for a recent review paper on applications of AAA.

Consider $N>0$ distinct sampling points in the set $\cT = \{\tau_1, \tau_2, \ldots,\tau_N\}$ together with $N$ function evaluations at these points, denoted with $\{f_1, f_2, \ldots,f_N\}$.
The AAA algorithm is based on an iteration. 
Let the order of the rational approximant $\raaa(\varIC)$ that is constructed through the algorithm in \cite{nst18} at step $\ell \geq 1$ be $(\ell,\ell)$ (degree of numerator and denominator), which is written in the following barycentric format: 
\begin{equation}\label{eq:AAA_approx}
\raaa(\varIC) = \frac{ \sum_{k=0}^\ell \frac{ \alpha_k w_k }{\varIC - \lambda_k }}{ \sum_{k=0}^\ell \frac{\alpha_k }{\varIC- \lambda_k}},
\end{equation}
where $\{\lambda_0,\ldots,\lambda_\ell\}$ are the interpolation points (selected from the set $\cT$), $\{w_0,\ldots,w_\ell\}$ are the corresponding function evaluations, and $\{\alpha_0,\ldots,\alpha_\ell\}$ are the weights. The method enforces interpolation at the $\ell$ points $\{\lambda_0,\ldots,\lambda_\ell\}$, i.e. $\raaa(\lambda_k) = w_k = F(\lambda_k), \, \forall \, 0 \leqslant k \leqslant \ell$. Then, at each step, one needs to determine suitable weights $\alpha_0,\alpha_2,\ldots,\alpha_\ell$ by formulating an approximation problem. Instead of solving the more general problem above (which is nonlinear in variables $\alpha_1, \alpha_2, \ldots, \alpha_\ell$), one solves a relaxed problem by means of linearization. This procedure leads to a least squares minimization problem that involves a (rectangular) Loewner matrix of dimension $(N-\ell) \times \ell$. The solution relies on computing the SVD of the latter, and selecting the weights as the entries of the (normalized) last right singular vector.
At the end of the iteration, the rational approximant is of order $(r,r)$; one can also explicitly write the AAA function as $\raaa(\varIC) = \hat{\bC}_{\mathrm{AAA}} (\varIC \hat{\bE}_{\mathrm{AAA}}-\hat{\bA}_{\mathrm{AAA}})^{-1} \hat{\bB}_{\mathrm{AAA}}$,
where the matrices are of dimension $r=\ell+1$, given explicitly in terms of the selected support points  $\{\lambda_0,\ldots,\lambda_r\}$, the function evaluations $\{w_0,\ldots,w_r\}$ and the weights $\{\alpha_0,\ldots,\alpha_r\}$. Note that for the numerical examples presented in what follows, we will be using the AAA implementation available in the \chebfun~ tool \cite{chebfun}\footnote{Also available at \url{https://github.com/chebfun/chebfun}.}.

Computing rational minimax approximations is a challenging task, especially when there are singularities on or near the interval of approximation. In the recent contribution \cite{fntb18}, the authors propose different combinations of robust algorithms that are based on rational barycentric representations. One key feature is that the support points are chosen adaptively as the approximant is computed. The \chebfun~ minimax implementation is available in  \cite{chebfun}. The  AAA-Lawson algorithm in \cite{nakatsukasa2020algorithm} uses an iteratively re-weighted least-squares method followed by a classical Remez algorithm to enforce  quadratic convergence. Minimax approximants can be recognized by their equioscillatory error curves in real cases and nearly circular error curves in complex problems.


In the recent note \cite{trefethen2024computation}, the authors present an algorithm based on AAA-Lawson approximation~\cite{nakatsukasa2020algorithm}, enhanced by various modifications that make it suitable for dealing with functions of the flavor of $f(\varIC) = \hbox{sign}(\varIC)$ and for ensuring convergence of the Lawson phase. This variant/improvement of the standard AAA algorithm in \cite{nst18} is meant as a reliable tool for numerical computations of solutions for the third and fourth Zolotarev problems, i.e., classical problems in approximation theory.


An attentive reader may realize that in \cite{trefethen2024computation}, these developments are deemed a "breakthrough" since "a reliable method for computing these functions has not been available before". However, the LF is known since 1986 in its one-sided version \cite{aa86} and since 2007 in its two-sided version \cite{ajm07} (the latter was used in this work).


\section{The two symmetric circles case}
\label{sec:1a}

\subsection{Experimental setup and results} 

We start the numerical exposition with a rather simplified configuration, namely the two symmetric circles case. This configuration is also denoted \texttt{'1a'} in \cite{trefethen2024computation}. For this topology, the true optimal rational approximation solution for the Zolotarev problem Z4 is explicitly known. Let us consider two circles $E$ ($-1$) and $F$ ($+1$) symmetric (with respect to the imaginary axis), centered along the horizontal axis, with radius $\rho$ and centered at $\pm \alpha$. In this case, the optimal rational $r$-th order approximant $r_r^*(\varIC)$ derived in \cite{starke1992near} is given by:
\begin{equation}
r_r^*(\varIC) = \left(\dfrac{\varIC-\sqrt{\alpha^2-\rho^2}}{\varIC+\sqrt{\alpha^2-\rho^2}}\right)^r
\text{ and }
\sigma_r = \left(\dfrac{1-\sqrt{\alpha^2-\rho^2}}{1+\sqrt{\alpha^2-\rho^2}}\right)^r.
\label{eq:case_1a_opt}
\end{equation}

The case considered in what follows corresponds to the one in \cite[Figure 1(a)]{trefethen2024computation},  for which $\rho = 1/2$ and $\alpha = 1$ (see also \cref{fig:intro:1a} and \cref{fig:intro:1a_1e14}). To compare the rational approximation obtained for problem Z4 using either LF or AAA, the following two scenario will be used (see also \cref{sec:num} for LF code description):
\begin{itemize}
\item[S1] (i) the LF, (ii) the AAA, (iii) the AAA-sign(damping) and (iv) the AAA-sign(damping)-Lawson (with 200 iterations)\footnote{The AAA results are obtained with optional arguments  \texttt{'sign',0,'lawson',0}, \texttt{'sign',1,'damping',.95} and \texttt{'sign',1,'damping',.95,'lawson',200}, using the AAA code in \chebfun.}. 
\item[S2] (i) the LF, (ii) the AAA-Lawson (with 200 iterations) and (iii) the AAA-Lawson (with 1000 iterations)\footnote{The AAA-Lawson is obtained with optional arguments  \texttt{'sign',1,'damping',.95,'lawson',200} and \texttt{'sign',1,'damping',.95,'lawson',1000}, using the AAA code in \chebfun.}. 
\end{itemize}

\def\nLocal{4} 
\subsection{Structure analysis of the rational solutions}

\begin{figure}[ht!]
    \centering
    \includegraphics[width=.75\textwidth]{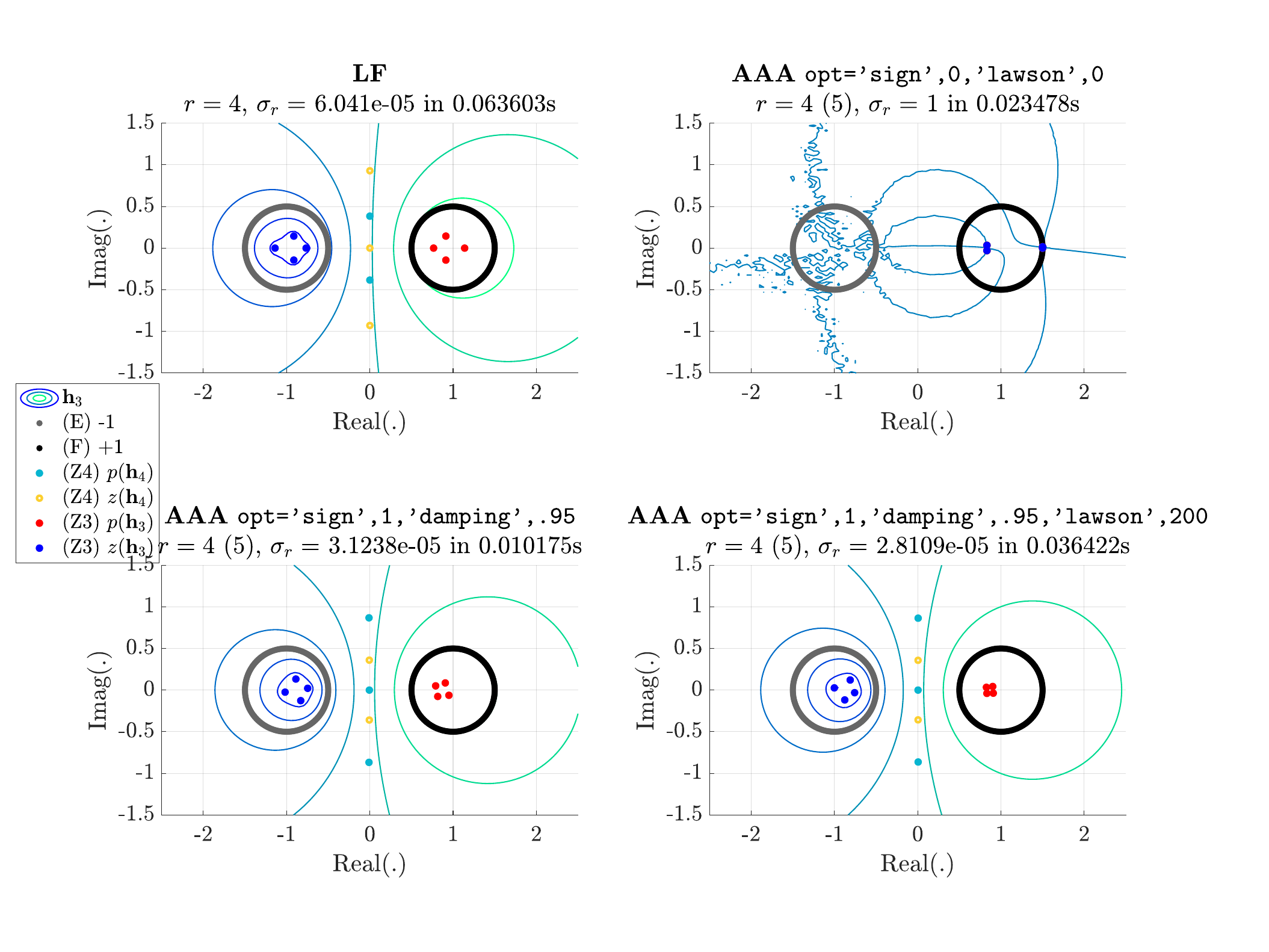}
    \caption{Case \texttt{'1a'}, with $r=\nLocal$: different approximation methods: LF (top left), AAA (top right), AAA with damping (bottom left), and AAA-Lawson with damping and  200 iterations (bottom right). Each frame shows: the Zolotarev rational functions of degree $r$ for Z3 ($\bh_3$) and Z4 ($\bh_4$) defined by connected sets $E$ (on the left) and $F$ (on the right); the Z3 function approximate contours levels $\log_{10} |\bh_3(\varIC)|= 0,-1,-2,\cdots,-30$ between the two domains; yellow circles and cyan blue dots respectively mark zeros and poles of the sign function $\bh_4$; blue and red dots respectively denote the zeros and poles of $\bh_3$. The minimum values $\sigma_r$ and computation time are given in the title, as well as the options of each AAA method. The number of singularities targeted (and obtained).}
    \label{fig:intro:1a}
\end{figure}

We begin with a simple approximation with $r=\nLocal$. The \cref{fig:intro:1a} shows the topology, approximation and poles/zeros distributions for the scenario S1. More specifically,  we display the rational approximation for Z3 with contours, and both the poles and the zeros for Z3 and Z4. Each frame illustrates the rational approximation reached by each method.  Comparing solely the approximation metric $\sigma_\nLocal$ (in each frame titles), this figure shows that every methods but AAA (with no option) is able to achieve a good approximation. In addition, it also clearly emphasizes the superiority of AAA combined with sign and Lawson iterations, both  being very close to the optimal one given in  \cref{eq:case_1a_opt}.

Now, let us focus on the AAA parametrizations leading to a good approximation and consider scenario S2. In the table below, we report the coefficients of both numerator and denominator (both real and imaginary parts) of the rational approximation obtained with each different method, and compare to the optimal Z4 approximant of the same order, obtained with \cref{eq:case_1a_opt} ($r=\nLocal$)\footnote{The caption mentions \texttt{cpv MACA64 MR2023b} since these results have been obtained by C. Poussot-Vassal on a Mac OS with \texttt{MATLAB} 2023b. See \cref{app:1a} for other architectures.}. The maximum approximation errors $\sigma_\nLocal$ are displayed for all approximants. 
\begin{table}[H] \tiny $$ 
\begin{array}{c||c>{\columncolor{colhl}}c|c>{\columncolor{colhl}}c|c>{\columncolor{colhl}}c|c>{\columncolor{colhl}}c|c>{\columncolor{colhl}}c|c>{\columncolor{colhl}}c|}  & \multicolumn{2}{c}{\text{Opt., $\sigma_{4}=2.7 \cdot 10^{-05}$}}  & \multicolumn{2}{c}{\text{LF, $\sigma_{4}=6 \cdot 10^{-05}$}}  & \multicolumn{2}{c}{\text{AAA-sign-L200, $\sigma_{4} =2.8 \cdot 10^{-05}$}}  & \multicolumn{2}{c}{\text{AAA-sign-L1000, $\sigma_{4} =2.8 \cdot 10^{-05}$}}   \\ \hline \text{Basis} & \text{real} & \text{imag.}  & \text{real} & \text{imag.}  & \text{real} & \text{imag.}  & \text{real} & \text{imag.}  \\   \hline  \hline z^4 & 0 & 0 & 0 & 0 & -140.0 & 16.0 & -150.0 & 7.3\\ z^3 & 3.5 & 0 & 3.7 & 3.3 \cdot 10^{-16} & 4.7 & -0.29 & 4.8 & -0.22\\ z^2 & 0 & 0 & -4.6 \cdot 10^{-15} & 3.3 \cdot 10^{-15} & -650.0 & 70.0 & -660.0 & 33.0\\ z & 2.6 & 0 & 3.2 & 5.6 \cdot 10^{-16} & 5.5 & -0.52 & 5.6 & -0.34\\ 1.0 & 0 & 0 & -1.3 \cdot 10^{-15} & -1.9 \cdot 10^{-15} & -80.0 & 8.8 & -82.0 & 4.1\\  \hline  \hline z^4 & 1.0 & 0 & 1.0 & 0 & 1.0 & 0 & 1.0 & 0\\ z^3 & 0 & 0 & -2.7 \cdot 10^{-15} & -1.3 \cdot 10^{-15} & -500.0 & 54.0 & -510.0 & 25.0\\ z^2 & 4.5 & 0 & 5.2 & 4.8 \cdot 10^{-16} & 7.9 & -0.7 & 7.9 & -0.49\\ z & 0 & 0 & -3.5 \cdot 10^{-15} & -2.2 \cdot 10^{-15} & -370.0 & 41.0 & -380.0 & 19.0\\ 1.0 & 0.56 & 0 & 0.75 & 9.5 \cdot 10^{-17} & 1.4 & -0.13 & 1.4 & -0.078  \\ \hline \end{array} 
 $$\normalsize \caption{Case \texttt{1a} (\texttt{cpv MACA64 MR2023b}), $r=4$, Z4: numerator (first lines) and denominator (last lines) coefficients.} \label{tab:sym-1a-r4-cpv_MACA64_MR2023b} \end{table}

The first comment for this case is that the optimal rational approximation \cref{eq:case_1a_opt} is attained with a real-valued function and a rational complexity of degree $(r-1,r)$. Therefore, for $r=\nLocal$, one notices that:
\begin{itemize}
    \item LF provides a solution with real-valued coefficients only (up to machine precision), while both AAA-sign-Lawson end up with non-negligible imaginary coefficients;
    \item by inspecting the numerator (first row), both AAA-Lawson suggest a non-null coefficient along with the $\varIC^\nLocal$ basis entry, resulting in a rational complexity of $R_{\nLocal,\nLocal}$ while the LF indicates a null term, leading to a complexity $R_{\fpeval{\nLocal-1},\nLocal}$ (as the optimal one);
    \item by inspecting the numerator along with the even basis entries $\varIC^0,\varIC^2,\dots$ and denominator of along with odd  basis entries $\varIC^1,\varIC^3,\dots$, it appears that the optimal have null coefficients. Again, only LF recovers this property (up to machine precision) while AAA shows high magnitude coefficients.
    \item by inspecting the $\sigma_\nLocal$ values for each method, AAA-Lawson provides a better approximant (i.e., with a lower $\sigma_\nLocal$) than LF (which is better than standard AAA); moreover, AAA-Lawson actually provides the near-optimal $\sigma_\nLocal$.
\end{itemize}
Choosing different order $r$, one may observe that the previous remarks still hold true\footnote{We refer to the \url{https://github.com/cpoussot/zolotarev} page for an in-depth analysis.} (see extensive report in \cref{app:1a}). Continuing the analysis, the following tables report the zeros and poles associated to the numerator and denominator for each rational approximation. 
\begin{table}[H] \tiny $$ 
\begin{array}{c||c|c|c|c|>{\columncolor{colhl}}c}  & \multicolumn{1}{c}{\text{Opt., $\sigma_{4}=2.7 \cdot 10^{-05}$}}  & \multicolumn{1}{c}{\text{LF, $\sigma_{4}=6 \cdot 10^{-05}$}}  & \multicolumn{1}{c}{\text{AAA-sign-L200, $\sigma_{4} =2.8 \cdot 10^{-05}$}}  & \multicolumn{1}{c}{\text{AAA-sign-L1000, $\sigma_{4} =2.8 \cdot 10^{-05}$}}   \\ \hline  1.0 & \mathrm{NaN} & \mathrm{NaN} & 0.00403-0.358{}\imath & 0.00404+0.358{}\imath\\ 2.0 & 0 & 4.1 \cdot 10^{-16}+5.92 \cdot 10^{-16}{}\imath & 0.00406+0.358{}\imath & 0.00402-0.358{}\imath\\ 3.0 & 0.866{}\imath & 3.69 \cdot 10^{-16}+0.929{}\imath & 0.0128-2.08{}\imath & 0.0127-2.09{}\imath\\ 4.0 & -0.866{}\imath & 4.46 \cdot 10^{-16}-0.929{}\imath & 0.0117+2.09{}\imath & 0.0116+2.09{}\imath  \\ \hline \end{array}
 $$\normalsize \caption{Case \texttt{1a} (\texttt{cpv MACA64 MR2023b}), $r=4$, Z4: zeros.} \label{tab:zer-1a-r4-cpv_MACA64_MR2023b} \end{table}
\begin{table}[H] \tiny $$ 
\begin{array}{c||c|c|c|c|>{\columncolor{colhl}}c}  & \multicolumn{1}{c}{\text{Opt., $\sigma_{4}=2.7 \cdot 10^{-05}$}}  & \multicolumn{1}{c}{\text{LF, $\sigma_{4}=6 \cdot 10^{-05}$}}  & \multicolumn{1}{c}{\text{AAA-sign-L200, $\sigma_{4} =2.8 \cdot 10^{-05}$}}  & \multicolumn{1}{c}{\text{AAA-sign-L1000, $\sigma_{4} =2.8 \cdot 10^{-05}$}}   \\ \hline  1.0 & -1.94 \cdot 10^{-16}+2.09{}\imath & 3.21 \cdot 10^{-16}-0.384{}\imath & 0.00377+6.65 \cdot 10^{-5}{}\imath & 0.00375-2.07 \cdot 10^{-5}{}\imath\\ 2.0 & -1.94 \cdot 10^{-16}-2.09{}\imath & 3.09 \cdot 10^{-16}+0.384{}\imath & 0.00535-0.864{}\imath & 0.00513+0.864{}\imath\\ 3.0 & -1.39 \cdot 10^{-17}+0.359{}\imath & 1.12 \cdot 10^{-15}-2.25{}\imath & 0.00516+0.864{}\imath & 0.00532-0.864{}\imath\\ 4.0 & -1.39 \cdot 10^{-17}-0.359{}\imath & 9.09 \cdot 10^{-16}+2.25{}\imath & 499.0-54.4{}\imath & 507.0-25.2{}\imath  \\ \hline \end{array}
 $$\normalsize \caption{Case \texttt{1a} (\texttt{cpv MACA64 MR2023b}), $r=4$, Z4: poles.} \label{tab:pol-1a-r4-cpv_MACA64_MR2023b} \end{table}

These tables emphasize that both the poles and zeros obtained with AAA-Lawson are quite far from the optimal ones and contain a non-negligible real part, while the LF ones are very close. In addition, as the AAA rational approximation coefficients have large imaginary parts, poles and zeros do not come in complex conjugate form, while LF results in almost perfect closed form. Then, one pole comes with a very large magnitude, far from the domain of interest. Finally, as AAA constructs a rational approximation of the type $R_{\nLocal,\nLocal}$, an additional zero appears.

\subsection{Accuracy and computational effort}

\begin{figure}[ht!]
    \centering
    \includegraphics[width=.75\textwidth]{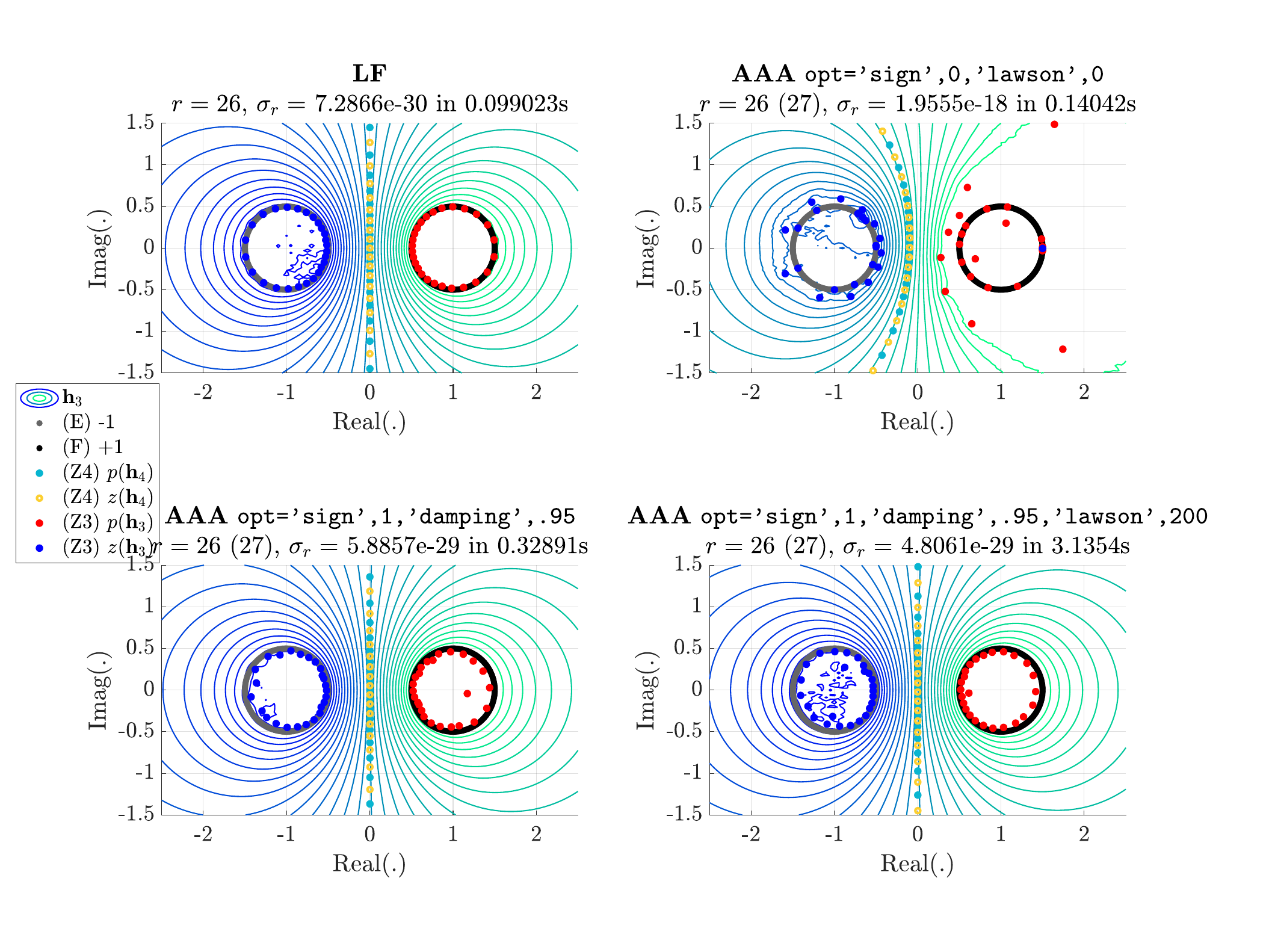}
    \caption{Case \texttt{'1a'}, with $r=26$: same as \cref{fig:intro:1a}.}
    \label{fig:intro:1a_1e14}
\end{figure}

Let now continue with approximation for varying rational degrees. First we compare the approximation features for significantly higher orders, using again scenario S1. In \cref{fig:intro:1a_1e14}, the approximation order $r$ has been selected such that the normalized singular values of the LF are below $10^{-14}$ (as suggested in \cref{eq:loew-rank}); then, the same objective is used for the three other methods.

What is remarkable in \cref{fig:intro:1a_1e14}, for the LF case, is the neat distribution of the poles and zeros of both Z3 and Z4 problems along the $\pm 1$ circles and vertical axis, with a perfect symmetry. This feature is guaranteed by the almost perfect alternate null and non-null coefficients and the null (machine precision) imaginary part. All other methods provide a slightly irregular distribution with imaginary part in the monomial coefficients. Moreover, in this setup, the optimal $\sigma_r\approx1.8144\cdot 10^{-30}$ \cref{eq:case_1a_opt} while the LF reaches the best approximation with $\sigma_n\approx 7.2866\cdot 10^{-30}$, very close to the optimal one. Notice also that the computation time is much lower for the LF with respect to the AAA versions, and even more with respect to the AAA-Lawson. 

\begin{figure}[ht!]
    \centering
    \includegraphics[width=.75\textwidth]{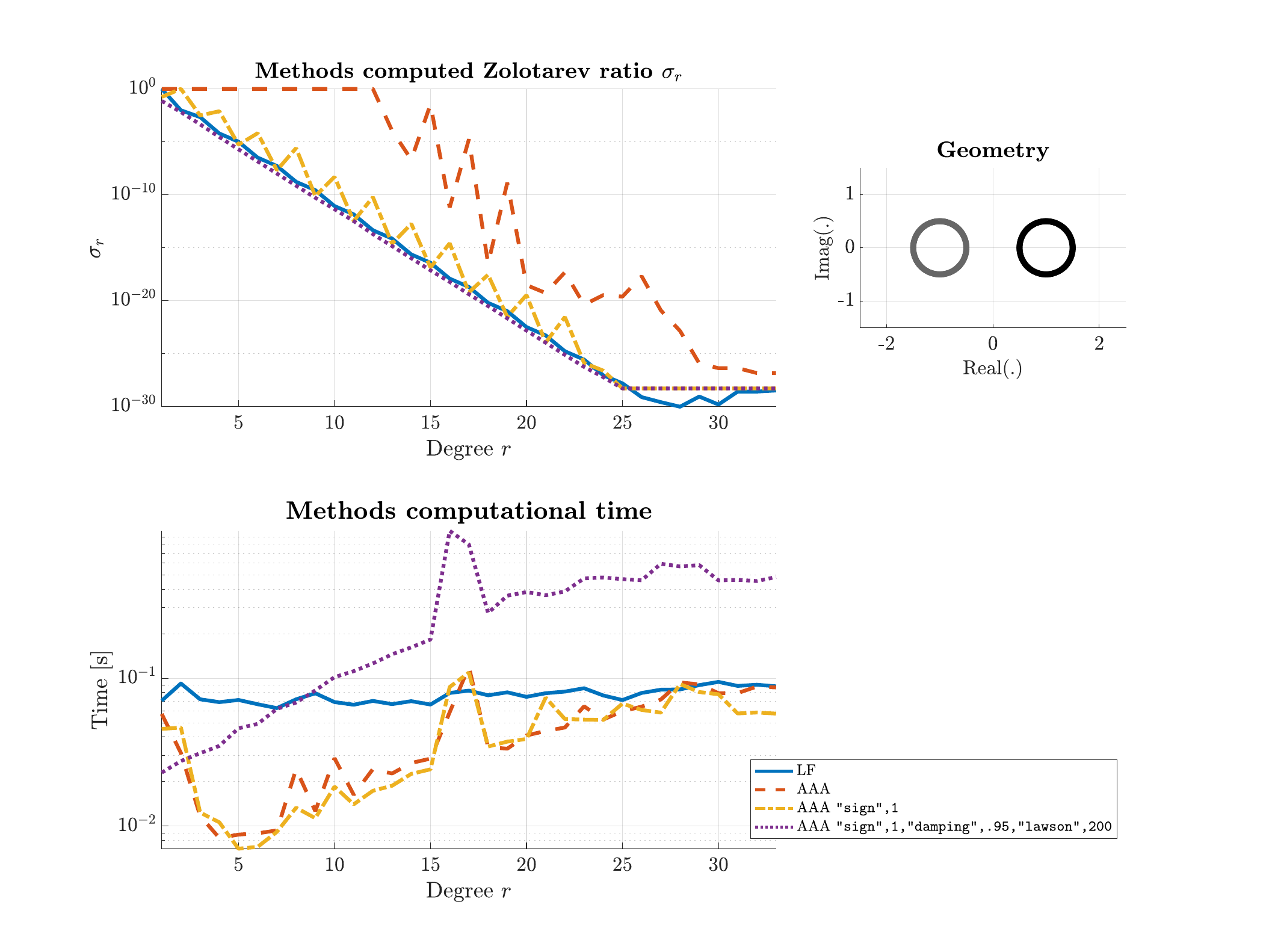}
    \caption{Case \texttt{'1a'}: computation time and accuracy for different approximation settings.}
    \label{fig:intro:1a_time}
\end{figure}

This computational time is indeed an important point from practical considerations. By now evaluating the Zolotarev ratio and computation time for different  rational approximation orders $r$, the \cref{fig:intro:1a_time} is obtained. One observes that the LF requires a constant computation time while the AAA variants are dependent on the approximation order $r$. Indeed, LF benefits from the fact that it is not iterative. 

\begin{remark}[Accelerating LF]
As the LF consists in (i) building the Loewner pencil \cref{eq:Loew_mat}, (ii) computing two SVDs \cref{eq:loew-rank}, and then (iii) truncating / projecting accordingly \cref{eq:Loew_red_mat}, in practice, steps (i) and (ii) may be computed once; then step (iii) must be computed for each truncation order. A direct consequence is the computational effort is limited to the first computation and the rest will be simply projections with a marginal cost. This procedure may not be applied to the AAA.
\end{remark}

Considering the accuracy, AAA-sign-Lawson generally provides a slightly better ratio number $\sigma_r$ for almost all $r$. However, it reaches a plateau for higher orders ($r>25$), which is anyway hard to distinguish with rounding numerical accuracy.

\subsection{Numerical robustness} 

\begin{figure}[ht!]
    \centering
    \includegraphics[width=.75\textwidth]{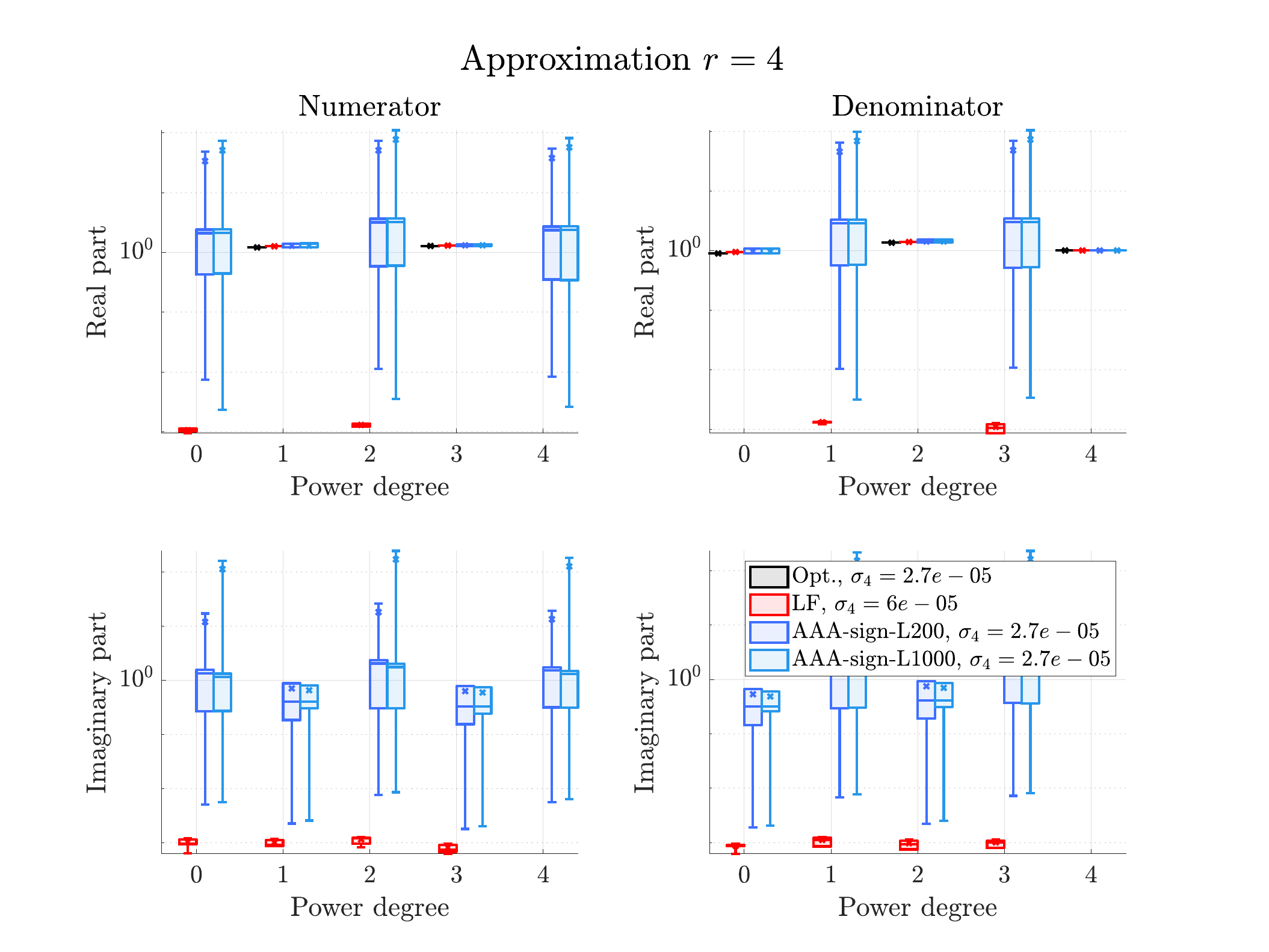}
    \caption{Case \texttt{1a}, $r=\nLocal$, Z4. Numerators and denominators real and imaginary coefficients dispersion over the \NARCH~ different platforms. For each frame, the $x$-axis refers to the monomial basis $\varIC$ degree, and the $y$-axis refers to the absolute value (dispersion) of the coefficient (in log-space). Left and right columns respectively show numerators and denominators, while upper and lower rows show the real and imaginary entries. Each bar shows statistics of each method: horizontal line are the maximal, minimal, median values; the cross is the mean value and shaded bar is the first quartile. Each color refer to a method.}
    \label{fig:intro:numden}
\end{figure}

We push our investigation by evaluating the results of the LF and AAA algorithms over different computer architectures settings, using scenario S2. More specifically, the code named \texttt{report\_1a\_step1}\footnote{Available at \url{https://github.com/cpoussot/zolotarev}.} has been executed on \NARCH~ different computers, with different architecture or \texttt{MATLAB} version (or both). Results were stored and analyzed via \texttt{report\_1a\_step2}, and dispersion of the numerator and denominator coefficients computed in \texttt{report\_1a\_step3}. This procedure applied to the case $r=\nLocal$, led to \cref{fig:intro:numden}. It shows the distribution of the rational approximation coefficients (numerator and denominator) according to the running machine. 

From \cref{fig:intro:numden} lower (imaginary) part, we immediately observe that in every case, the LF provides an almost null imaginary contribution, as indicated by the optimal solution (no bar), while AAA versions lead to large coefficients, with a spread distribution. From \cref{fig:intro:numden} upper (real) part, the alternate null/non-null coefficients is well restituted by LF, as the optimal solution, while AAA always provides non-null coefficients. Again, the AAA leads to a widely spread dispersion of all coefficients as a function of the architecture, which is puzzling.  

Still we point that in all cases, the approximation ration score $\sigma_r$ was constant and  good for all cases (S2 scenario), and close to the optimal.

\subsection{Preliminary conclusions and open questions} 

As highlighted in this simple case (of two circles), both AAA and AAA-Lawson variants produce {complex-valued approximants} and a complexity $R_{r,r}$, whereas the optimal result suggests a solution with real coefficients only and a complexity $R_{r-1,r}$. Alternatively, LF always produces  real-valued approximants (up to machine precision), without explicitly enforcing such a property\footnote{Notice that it is also possible to enforce parameters' realness, as explained in \cite[\S 2.4.5]{ali17}.}, and a complexity $R_{r-1,r}$. Second, only the LF seems to be able to recover the structure of the solution, i.e., the real coefficients and the alternate null and non-null coefficients in the numerator and denominator monomial coefficients, thus discovering the exact complexity and structure of the optimal approximant. Moreover, in AAA and AAA-Lawson, as $r$ increases, the spurious imaginary part of the coefficients becomes no longer negligible and even predominant. Finally, inspecting the quality of the approximation evaluated through the lens of $\sigma_r$, it appears that AAA is not accurate. Its conjunction with the Lawson iteration, however, seems to lead to significantly better results, i.e., close to the Zolotarev optimal, and better than those obtained with LF. The analysis of the poles and zeros of the Z4 approximations (and comparing with the optimal ones) showed, however, an inaccurate restitution of the AAA-Lawson, while the LF provides an accurate matching.

Along with these first comments, we believe it is legitimate to seek the significance in how spurious these complex terms are, and how they can provide a better $\sigma_r$ approximation. Indeed, while from a practical point of view, a better approximation (lower $\sigma_r$) may be satisfactory, from a theoretical point of view, improving accuracy by adding spurious coefficients is not.

Finally, we point that all experiments were run on \NARCH~ different computational environments (consisting of various operating systems and different versions of \texttt{MATLAB}). While the LF led to constant results, the AAA consistently led to significantly different results, both concerning the distribution of poles and zeros and the magnitude of the monomial coefficients appearing in the numerator and denominator. 

\section{Numerical experiments}\label{sec:num}

We now conduct numerical illustrations and a discussion of the effectiveness of the LF for solving the Z3 and Z4, and compare it to different AAA approaches through different (more complex) topologies\footnote{All the computations were performed on an Apple MacBook Air with 512 GB SSD and 16 GB RAM, with an M1 processor. The software used was MATLAB 2023b.} (more details and figures are available in \cref{app:all}). Let us first present a \texttt{MATLAB} code package\footnote{The \texttt{+zol} package is available at \url{https://github.com/cpoussot/zolotarev}.}. 

\subsection{Matlab code package, functions and code sample}

The Zolotarev code package embeds the examples reported in this work, the very basic (simplified) LF implementation for this setup, and a set of useful functions to approximate Z3 and Z4 problems. The package embeds the following functions that are sequentially used (see also the example next).

\begin{itemize}
\item \texttt{zol.example} is used to load the Zolotarev example: \\
\texttt{[pts,val,data] = zol.example(NAME)}\\
where \texttt{NAME} below  corresponds to a particular topology, and where output arguments contain the interpolation points (\texttt{pts}), corresponding values (\texttt{val}), and information regarding the example (\texttt{data}).
\begin{center}
\texttt{NAME=\{'1a', '1b', '1c', '1d', '1e', '1f', '2a', '2b', '2c', '2d', '3a', '3b', '3c', '3d', '7', 'spiral1', 'pm2'\}}    
\end{center}

\item \texttt{zol.example2data} is used to format the Zolotoarev (Z4) example to the LF:\\
\texttt{[la,mu,W,V] = zol.example2data(pts,val,data)}\\
where \texttt{pts,val,data} are the output arguments of \texttt{zol.example}. Then, output arguments contain the right  (\texttt{la,W}) and left (\texttt{mu,V}) interpolation points and values, being the main ingredients for the LF.

\item \texttt{zol.loewner} is the LF used for the Zolotarev (Z4) rational approximation.\\
\texttt{[h4,info] = zol.loewner(la,mu,W,V,opt)} \\
where, \texttt{la,mu,W,V} are the output arguments of \texttt{zol.example2data}. Then outputs gather the rational approximant of Z4 (handle function \texttt{h4}) and information about the LF (\texttt{info})\footnote{Gathers information about the LF process, as detailed in \cite{ali17}.}.
\item \texttt{zol.pb4\_to\_pb3} is used to convert the Zolotarev Z4 to Z3 problem:\\
\texttt{[h3,hp,hsig,htau] = zol.pb4\_to\_pb3(h4,pts,val)}\\
where \texttt{h4,pts,val} are output arguments of \texttt{zol.loewner} and \texttt{zol.example}. Then, outputs gather the rational approximant of Z3 (\texttt{h3}), the $p$, $\sigma$ and $\tau$ scalars (\texttt{hp,hsig,htau}), as defined in \cref{sec:zolo}.
\end{itemize}

We now present a brief code sample that solves the Z4 and Z3 problems for the two symmetric circle case (denoted \texttt{'1a'}). \\
\vspace{-2mm}
{\small
\begin{verbatim}
    % Load example '1a'
    [pts,val,data]  = zol.example('1a');
    % Select left/right interpolation points and values 
    [la,mu,W,V]     = zol.example2data(pts,val,data);
    % Apply basic LF to get h4, poles and zeros (Z4) 
    opt.target      = 1e-14; % order selected via singular value decay
    [h4,info]       = zol.loewner(la,mu,W,V,opt);
    h4poles         = eig(info.Ar,info.Er);
    h4zeros         = eig([info.Ar info.Br;info.Cr 0],blkdiag(info.Er,0));
    % Transform h4 (Z4) to h3 (Z3)
    [h3,hp,hsig]    = zol.pb4_to_pb3(h4,val,pts);
    h3poles         = eig([info.Ar info.Br;-info.Cr (hp)],blkdiag(info.Er,0));
    h3zeros         = eig([info.Ar info.Br; info.Cr (hp)],blkdiag(info.Er,0));
\end{verbatim}
}
\normalsize

\subsection{Numerical test cases}

We evaluate here the performance of multiple rational approximation methods over an extensive collection of Zolotarev topologies from \cite{trefethen2024computation} (more experiments are listed in \cref{app:all}). Here we concentrate in comparing the LF with AAA-sign(damping)-200 Lawson iterations\footnote{\texttt{'sign',1,'damping',.95,'lawson',200}.}.

\begin{figure}[ht!]
    \centering
    \includegraphics[width=.75\textwidth]{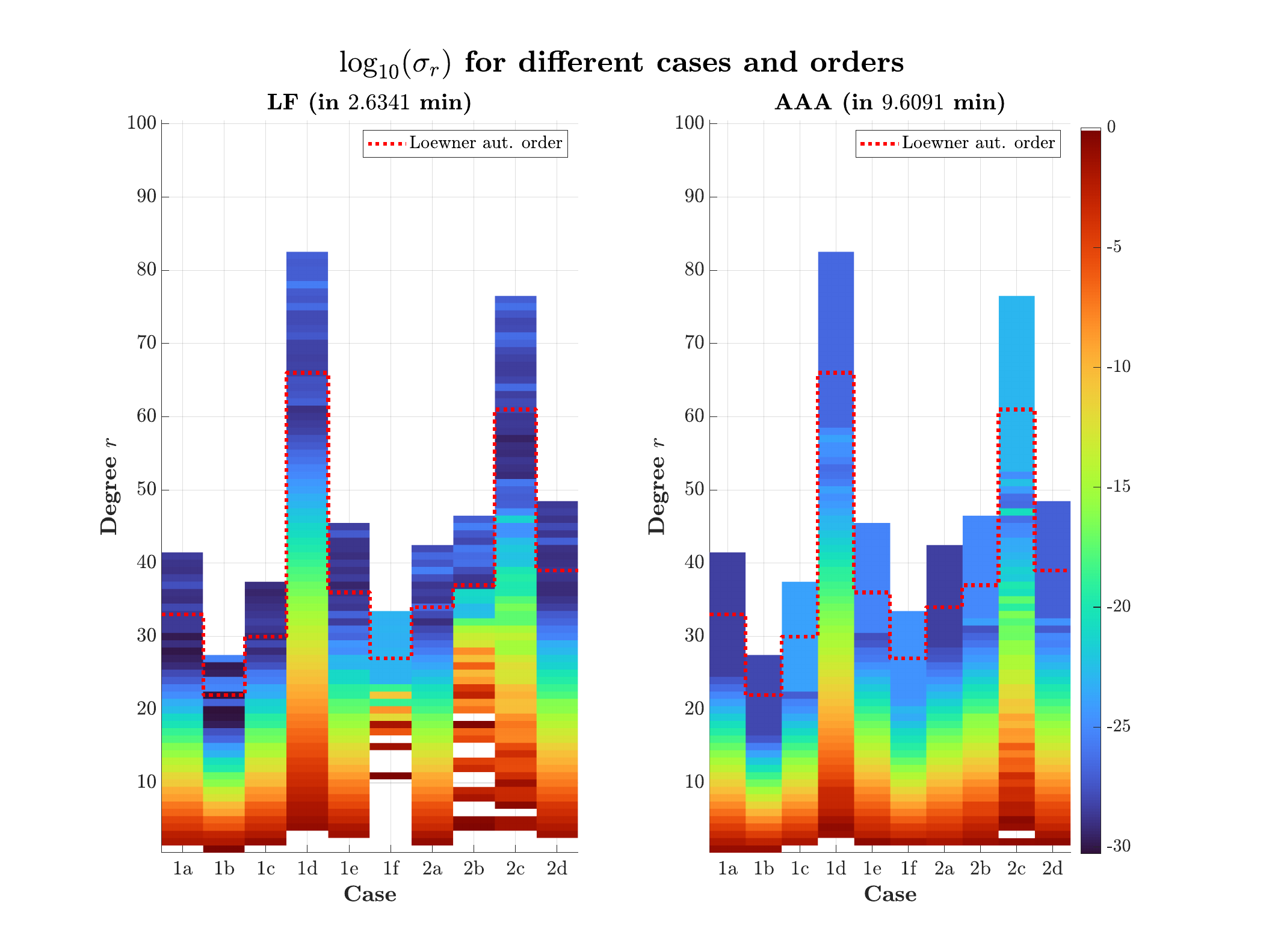}
    \caption{Z3 rational approximations for ten cases, with different approximation degrees $r$. Figure obtained by running \texttt{demo3\_LNT.m} file.}
    \label{fig:case_order}
\end{figure}

We conducted the following additional numerical experiments (accessible though the GitHub page \url{https://github.com/cpoussot/zolotarev} and the \texttt{demo3\_LNT.m} file), leading to  \cref{fig:case_order}. The latter shows the approximation $\sigma_r$ as a function of the approximation order ($y$-axis) for the ten first cases (groups 1 and 2). The dashed red lines correspond to the truncation order in the LF with a threshold fixed at $10^{-16}$, therefore, below machine precision. The color indicates the accuracy of the approximation in log-scale.

This figure shows two points. First, the LF (and AAA-Lawson) approximation perform well in approximating the Z4 problem. They do not fail for higher orders\footnote{Indeed, at this point, rounding errors are probably around, but the reason is not clear to us yet (as we chase for accuracy at $\sigma_r<10^{-20}$).}. Second, the following remarks can be made and would deserve additional  investigations: 
\begin{itemize}
    \item The overall table construction time is reported in the title showing that LF about three time faster than AAA-Lawson (with this setup, and options).
    \item Inspecting the \texttt{'1c'} case, AAA-Lawson behaves well up to 20 but then loses in accuracy, while LF provides multiple orders of magnitude better results.
    \item Inspecting the \texttt{'1f'} case, it is interesting to note that the LF fails for orders lower than about 20 while AAA-Lawson performs well and has a smooth and regular accuracy improvement with increasing orders. Moreover, for high orders, AAA-Lawson seems slightly more accurate. This case clearly favors of the AAA-Lawson strategy.
    \item Inspecting the \texttt{'2b'} case leads to a similar comment, but here, LF leads to better approximation when the order increases.
    \item In all other configurations, \cref{fig:case_order} shows that AAA-Lawson has difficulties when the order increases and get stuck to a constant approximation error while LF gets better (not monotonically) when reaching the singular value limit (at machine precision). Again here, rounding and \texttt{MATLAB} accuracy are around, yet to be understood.
\end{itemize}
More detailed and additional configurations are listed in \cref{app:all}.

\section{Conclusion}
\label{sec:conclusion}

In this paper, we have shown that the Loewner framework is a viable candidate for approximating the third and fourth Zolotarev problems. In particular, we have shown that the LF (i) is fast, accurate, and numerically robust (to the computational environment); (ii) approximates by compression with controlled complexity and accuracy, without either user intervention or iteration, and; (iii) is able to discover the polynomial and rational structure of the optimal solution. 

We have also provided a comprehensive comparison numerical study of LF with different versions of the AAA algorithm. In all configurations tested, the standard AAA variant was not able to accurately approximate the Z3 and Z4 problems. This pitfall was corrected thanks to the use of the sign constraint and the Lawson iteration (as in \cite{trefethen2024computation}). However, this was achieved at the price of a dedicated data treatment, leading to a notably higher computational cost, and an addition of poles and zeros and/or polynomial coefficients' imaginary parts, hard to interpret (and not in the correct optimal solution space, for the two symmetric circle case). Finally, we observed that the Lawson iteration is (i) not monotonically converging, which may result in a worse approximation than for AAA with sign only option (and without any iteration) and (ii) that results strongly vary with the computational environment. 

The code used for generating all results reported here (figures and tables) is made available at the following address:
\begin{center}
\url{https://github.com/cpoussot/zolotarev}
\end{center}

\subsection*{Acknowledgments}

We acknowledge A. dos Reis de Souza for running experiments.

\bibliographystyle{plain}
\bibliography{references}

\appendix

\section{Detailed results on the two symmetric circles case}\label{app:1a}
\foreach \n in {3,...,\CAS}{
\subsection{Approximation objective $r=\n$}
\subsubsection{Numerator and denominator coefficients}

The numerator and denominator coefficients (both real and complex parts) are given in the following tables. Each table gives the results obtained on different platform.
\tiny
\input{data/cpv_MACA64_MR2023b/tables/1a_nd_r\n}
\input{data/cpv_MACA64_MR2022b/tables/1a_nd_r\n}
\input{data/cpv_PCWIN64_MR2025b/tables/1a_nd_r\n}
\input{data/AlexReis_PCWIN64_MR2025b/tables/1a_nd_r\n}
\input{data/AlexReis_PCWIN64perso_MR2017a/tables/1a_nd_r\n}
\input{data/ivg_GLNXA64_MR2024a/tables/1a_nd_r\n}
\normalsize

\begin{figure}[H]
\centering
    \includegraphics[width=.75\textwidth]{figures/1a/num_den_r\n.pdf}
    \caption{Case \texttt{1a}, $r=\n$, Z4. Numerator and denominator dispersion over  the different platforms.}
    \label{fig:numden_\n}
\end{figure}

\subsubsection{Poles}

The poles are given in the following tables. Each table gives the results obtained on different platform.
\tiny
\input{data/cpv_MACA64_MR2023b/tables/1a_pol_r\n}
\input{data/cpv_MACA64_MR2022b/tables/1a_pol_r\n}
\input{data/cpv_PCWIN64_MR2025b/tables/1a_pol_r\n}
\input{data/AlexReis_PCWIN64_MR2025b/tables/1a_pol_r\n}
\input{data/AlexReis_PCWIN64perso_MR2017a/tables/1a_pol_r\n}
\input{data/ivg_GLNXA64_MR2024a/tables/1a_pol_r\n}
\normalsize

\subsubsection{Zeros}
The zeros are given in the following tables. Each table gives the results obtained on different platform.

\tiny
\input{data/cpv_MACA64_MR2023b/tables/1a_zer_r\n}
\input{data/cpv_MACA64_MR2022b/tables/1a_zer_r\n}
\input{data/cpv_PCWIN64_MR2025b/tables/1a_zer_r\n}
\input{data/AlexReis_PCWIN64_MR2025b/tables/1a_zer_r\n}
\input{data/AlexReis_PCWIN64perso_MR2017a/tables/1a_zer_r\n}
\input{data/ivg_GLNXA64_MR2024a/tables/1a_zer_r\n}
\normalsize


\subsubsection{Poles and zeros visualization (case \texttt{cpv MACA64 MR2023b})}

\begin{figure}[H]
\centering
    \includegraphics[width=.45\textwidth]{data/cpv_MACA64_MR2023b/figures/pz_r\n.pdf}
    \includegraphics[width=.45\textwidth]{data/cpv_MACA64_MR2023b/figures/pz_r\n_zoom.pdf}
    \caption{Case \texttt{1a} (\texttt{cpv MACA64 MR2023b}), $r=\n$, Z4. Left frame: zeros (left) and poles (right). Right frame: same but discarding the largest real entry.}
    \label{fig:pz_\n}
\end{figure}

\subsubsection{Approximation visualization}

\begin{figure}[H]
\centering
    \includegraphics[width=.75\textwidth]{data/cpv_MACA64_MR2023b/figures/eval_r\n.pdf}
    \caption{Case \texttt{1a} (\texttt{cpv MACA64 MR2023b}), $r=\n$, Z4. $E$ and $F$ error evaluations.}
    \label{fig:eval_\n}
\end{figure}
}

\section{Detailed results on multiple cases}\label{app:all}
We evaluate here the performance of multiple rational approximation methods over an extensive collection of Zolotarev topologies from \cite{trefethen2024computation}. Specifically, we compare the LF with AAA, AAA with sign, and AAA with sign, damping, and 200 Lawson iterations. In every case, the rational approximation order $r$ of the figure similar to \cref{fig:intro:1a} is selected by setting the lowest relative singular value threshold of the LF to $10^{-14}$, then the same objective is used for the other methods. Then, similarly to \cref{fig:intro:1a_time}, the Zolotarev ratio and computation time are reported for increasing orders $r$. In most cases, we note that the LF leads to an approximant in the form $R_{r-1,r}$ while AAA and its variants are in the form $R_{r,r}$.

\subsection{Case \texttt{'1b'}}

\begin{figure}[H]
    \centering
    \includegraphics[width=.75\textwidth]{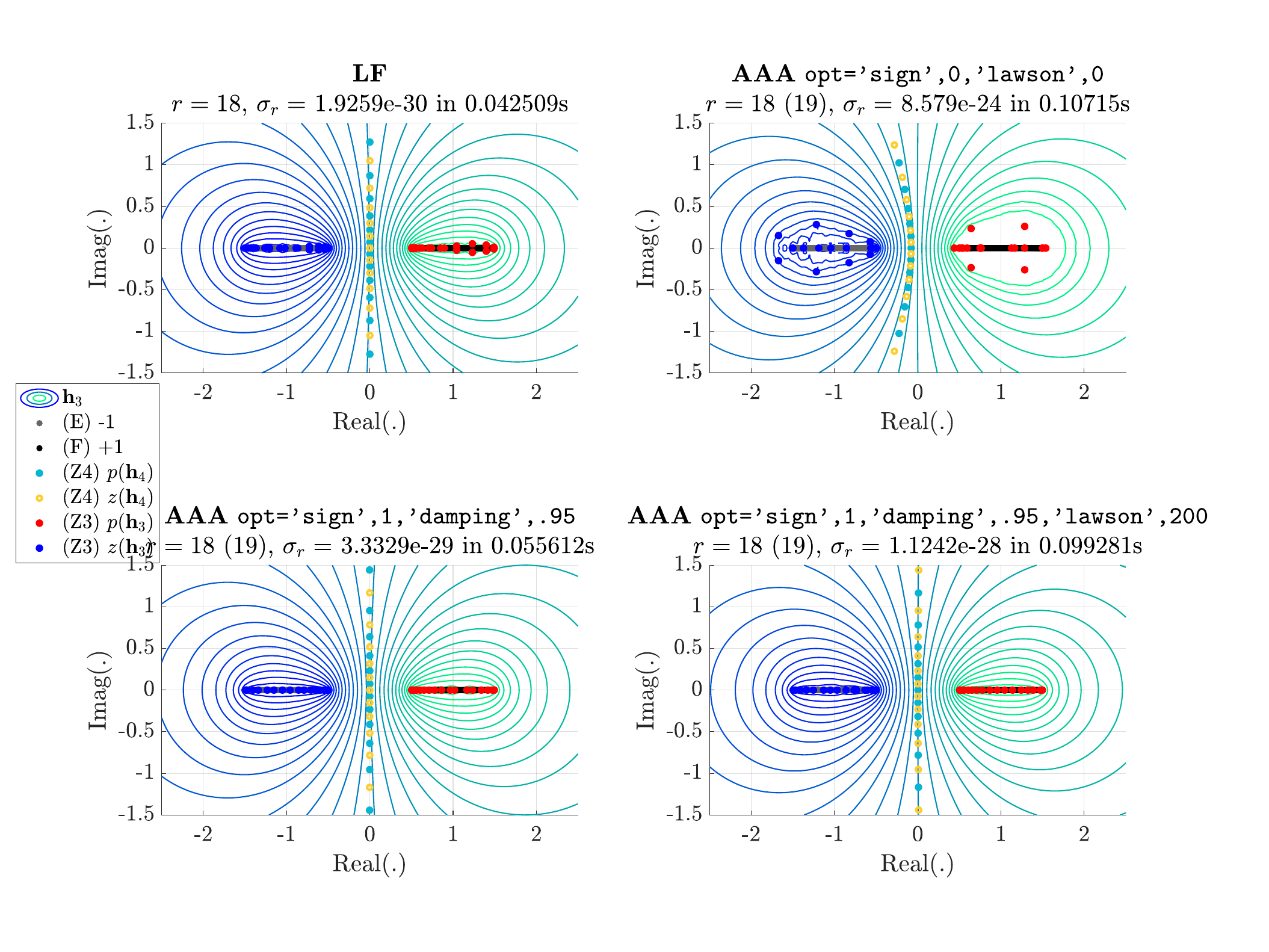}
    \includegraphics[width=.75\textwidth]{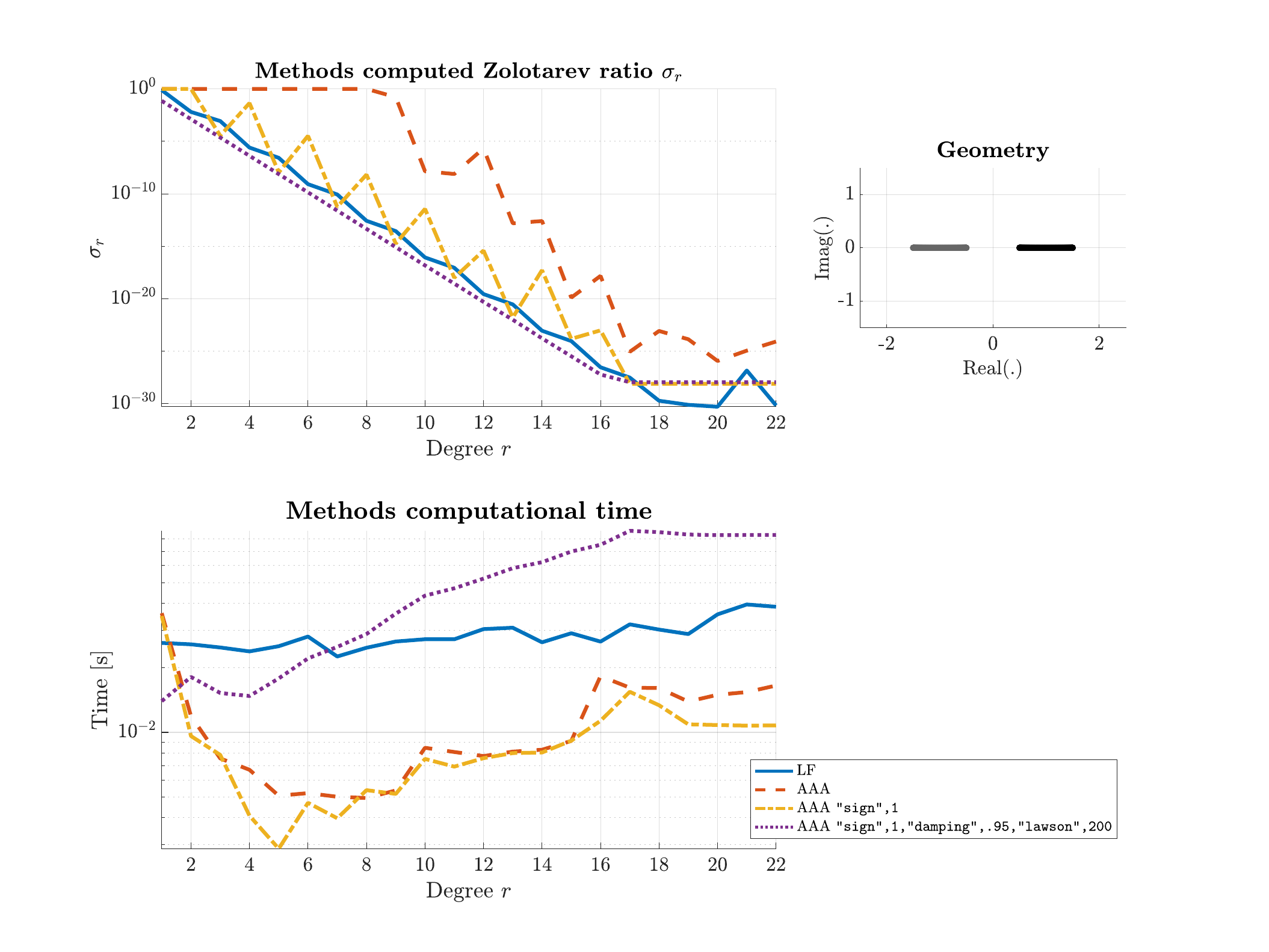}
    \caption{Case \texttt{'1b'}. Top: same description as \cref{fig:intro:1a}. Bottom: same description as \cref{fig:intro:1a_time}.}
    \label{fig:1b}
\end{figure}

\Cref{fig:1b} involves the classical sign function; it illustrates another symmetric case, i.e., a configuration where Z4 poles and zeros are aligned along the vertical axis, and Z3 poles and zeros are  necessarily real. In this case, the LF normalized singular values below the threshold indicate an approximation order of $r=18$. Then we use this order for each approximation method. By inspecting \cref{fig:1b}, the Zolotarev number indicates that AAA with sign and Lawson provides better results in almost all cases. LF performs better when orders are fairly high and also performs better than AAA and AAA signed. In addition, in this very simple case, the computation time of AAA and AAA with sign is below the LF. AAA with Lawson needs substantially more time.

\newpage
\subsection{Case \texttt{'1c'}}

\begin{figure}[H]
    \centering
    \includegraphics[width=.75\textwidth,align=c]{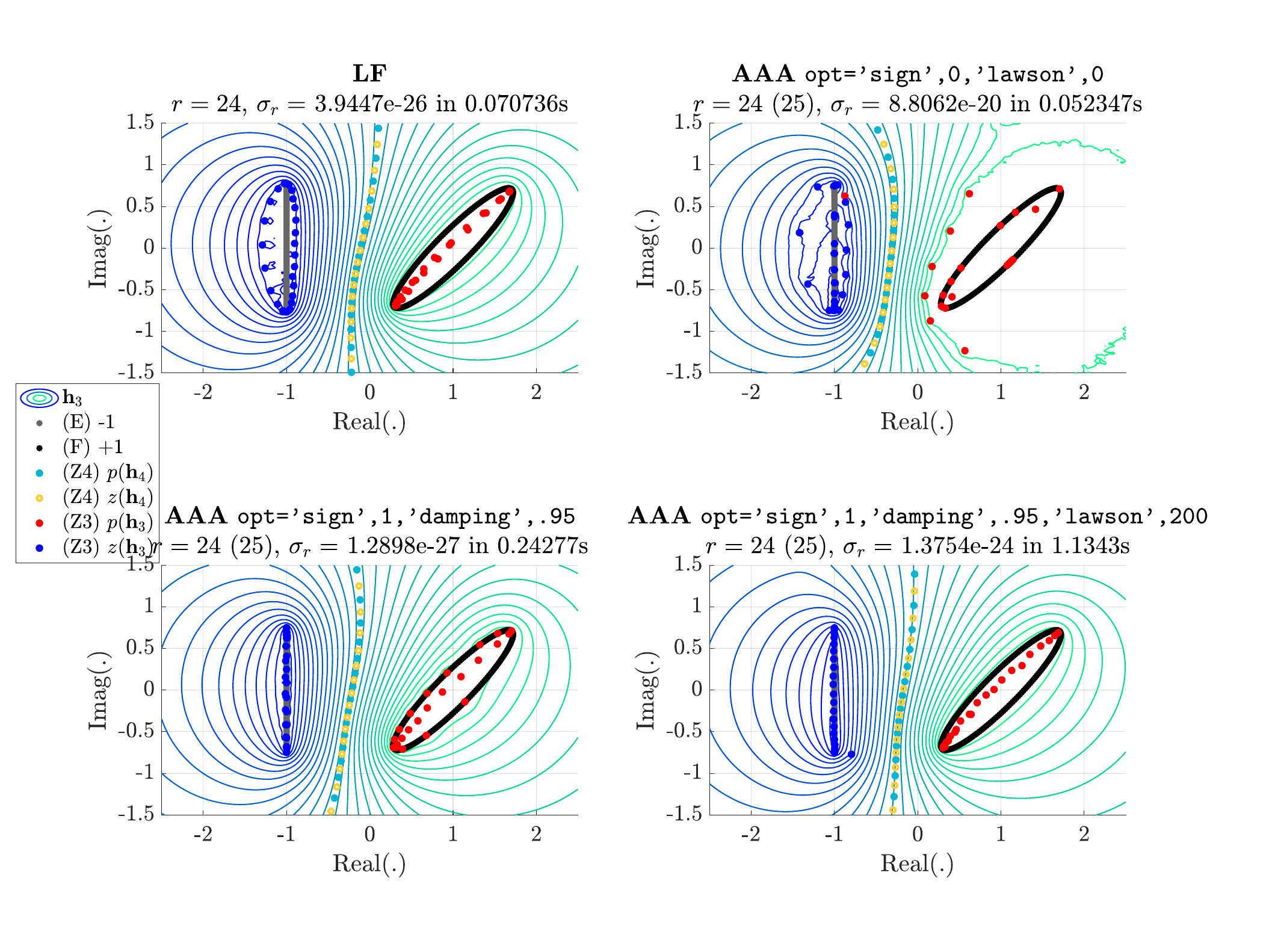}
    \includegraphics[width=.75\textwidth,align=c]{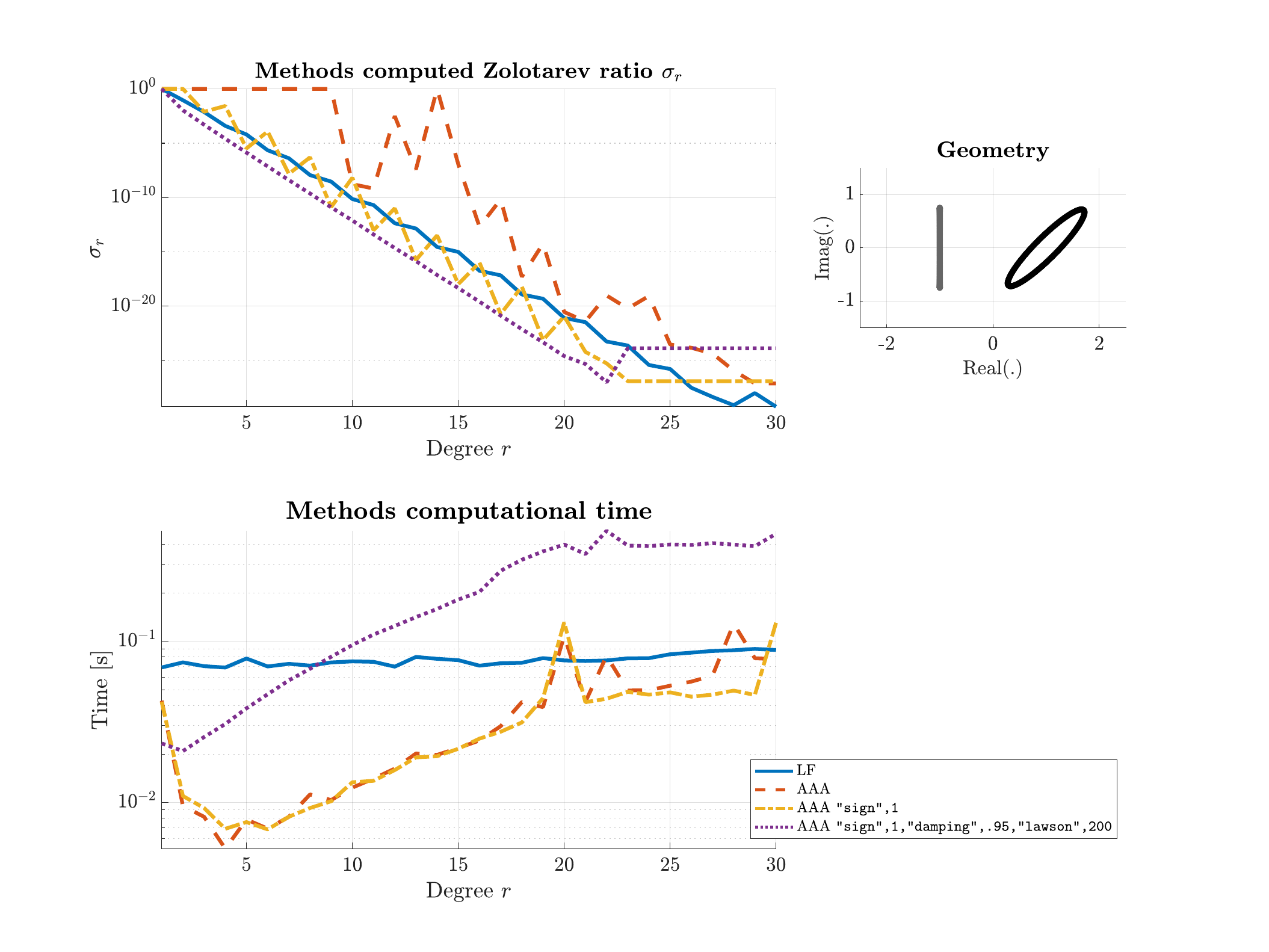}
    \caption{Case \texttt{'1c'}. Top: same description as \cref{fig:intro:1a}. Bottom: same description as \cref{fig:intro:1a_time}.}
    \label{fig:1c}
\end{figure}

\Cref{fig:1c} illustrates a non-symmetric case, i.e., a configuration where Z4 poles and zeros are not necessarily aligned along the vertical axis, and Z3 poles and zeros are not necessarily real and complex conjugated. In this case, the LF normalized singular values below the threshold indicate an approximation order of $r=24$. Then we apply this order to each approximation method. By inspecting  \cref{fig:1c} the Zolotarev number indicates that AAA with sign and AAA with sign and Lawson provides better results in almost all cases. LF outperforms the other methods when orders are fairly high.  In addition, in this simplified case, the computation time of AAA and AAA with sign is below the LF. AAA with Lawson needs substantially more time.

\newpage
\subsection{Case \texttt{'1d'}}

\begin{figure}[H]
    \centering
    \includegraphics[width=.75\textwidth,align=c]{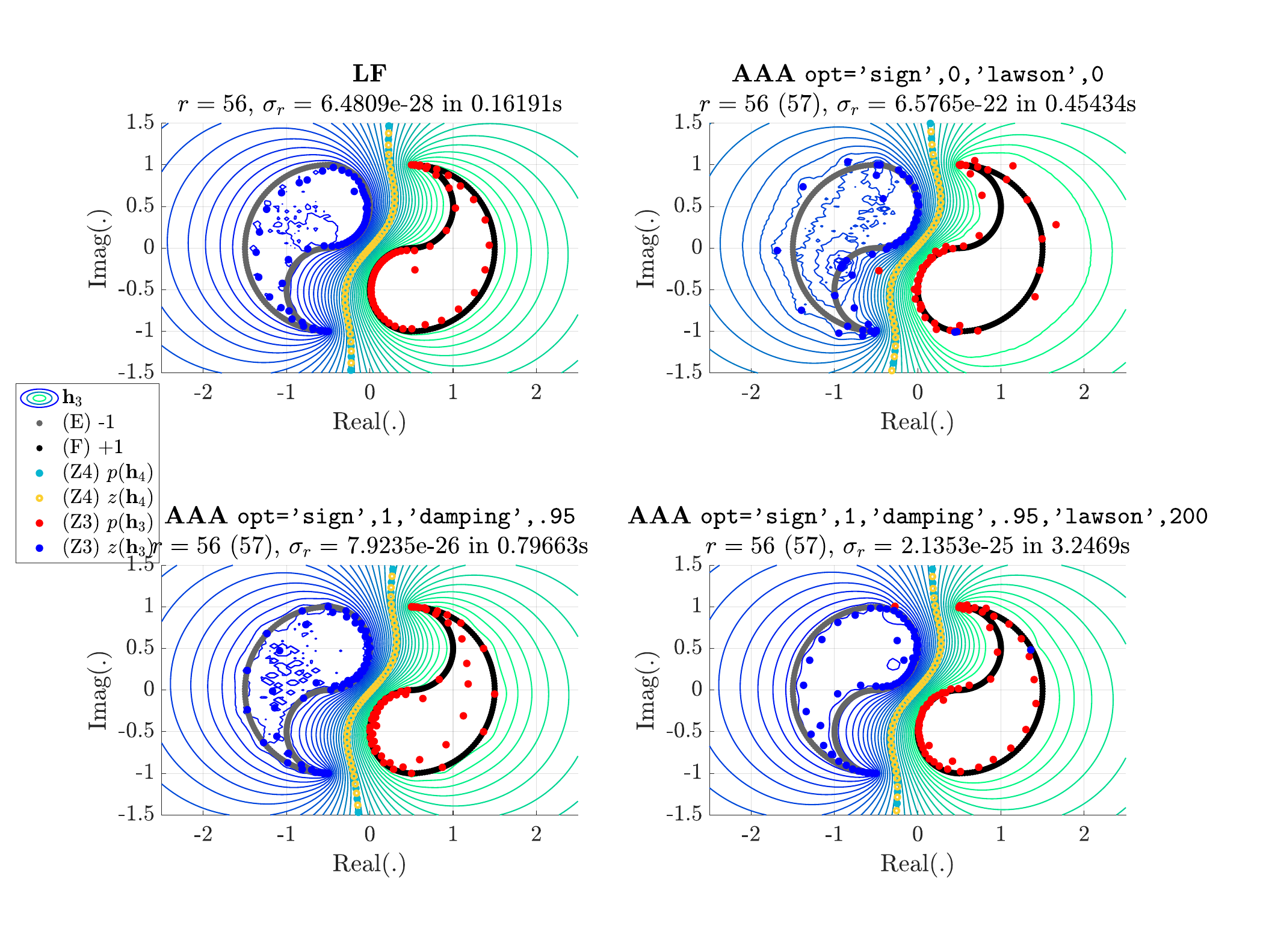}
    \includegraphics[width=.75\textwidth,align=c]{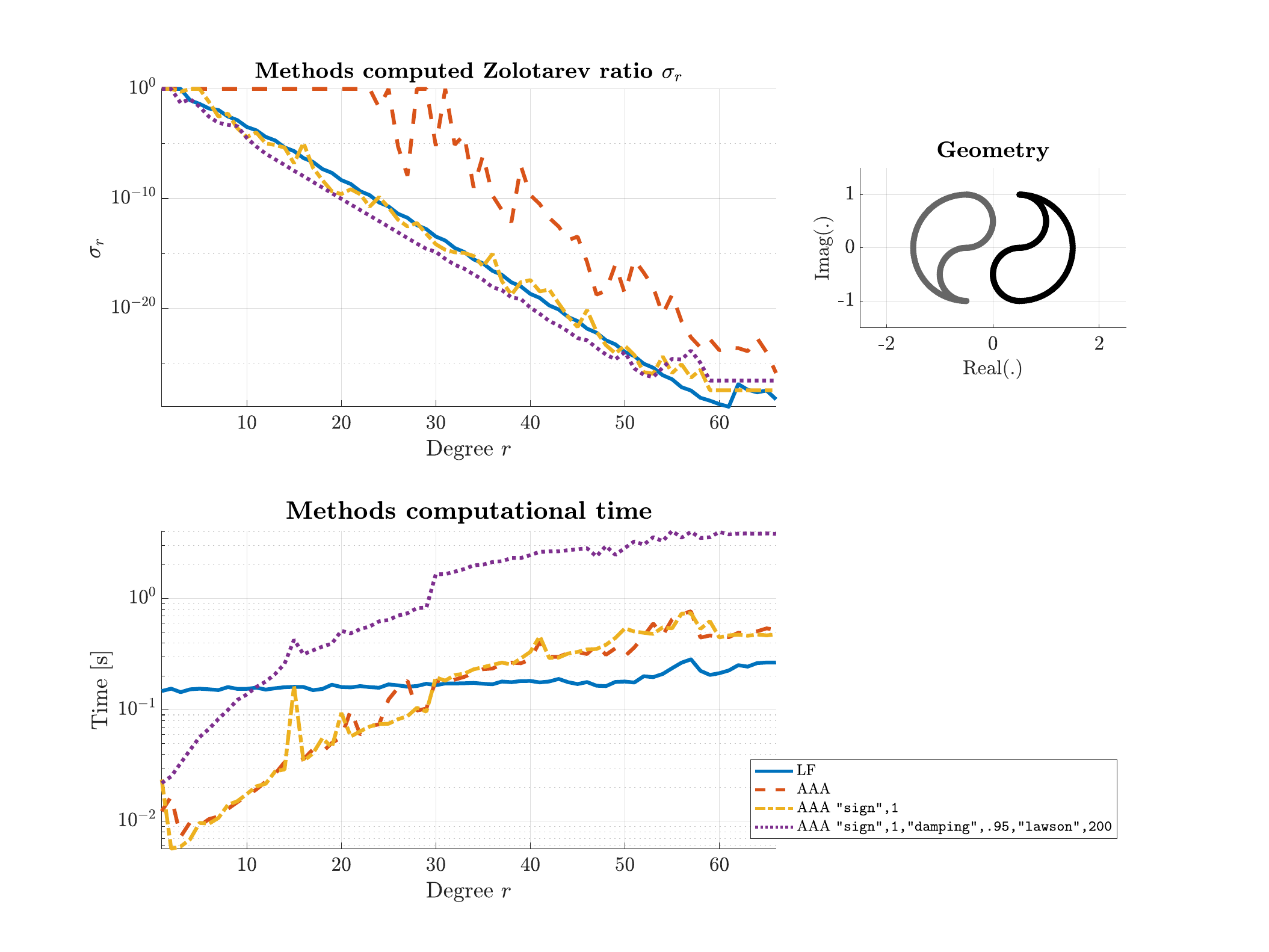}
    \caption{Case \texttt{'1d'}. Top: same description as \cref{fig:intro:1a}. Bottom: same description as \cref{fig:intro:1a_time}.}
    \label{fig:1d}
\end{figure}

This example considers a configuration symmetric by rotation. In \cref{fig:1d}, the time and accuracy plots indicate similar results in terms of Zolotarev number for LF, AAA sign, and AAA sign Lawson, with a slight advantage to AAA sign Lawson for almost all orders, except for high  approximation orders where LF outperforms. Moreover, inspecting the computation time rapidly shows an advantage of the LF. By inspecting the upper frame, only the LF and AAA sign result in poles and zeros totally separated and distributed along the yin/yang geometry. The two others show poles and zeros in the wrong part of the domain. 

\newpage
\subsection{Case \texttt{'2a'}}

\begin{figure}[H]
    \centering
    \includegraphics[width=.75\textwidth,align=c]{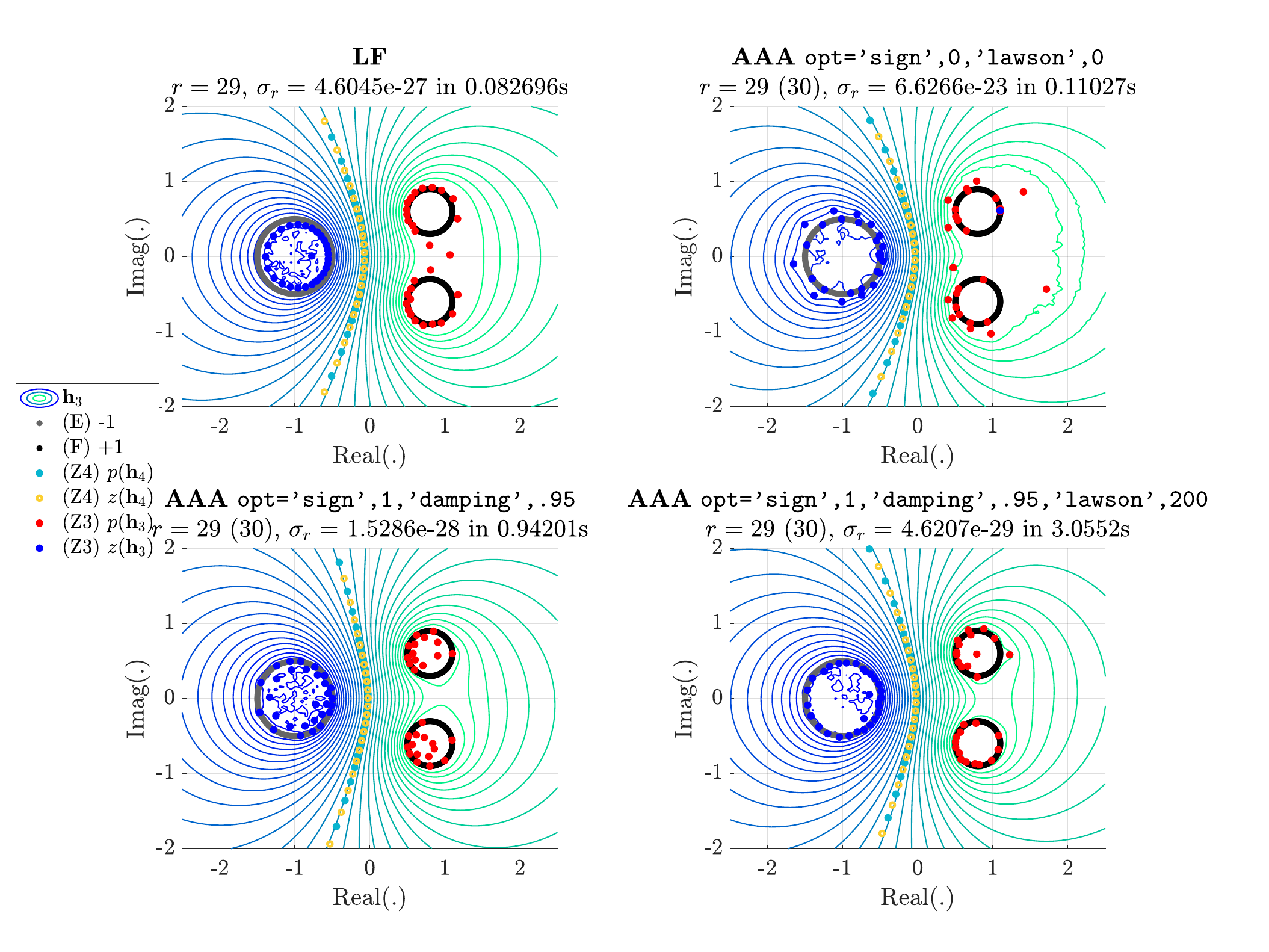}
    \includegraphics[width=.75\textwidth,align=c]{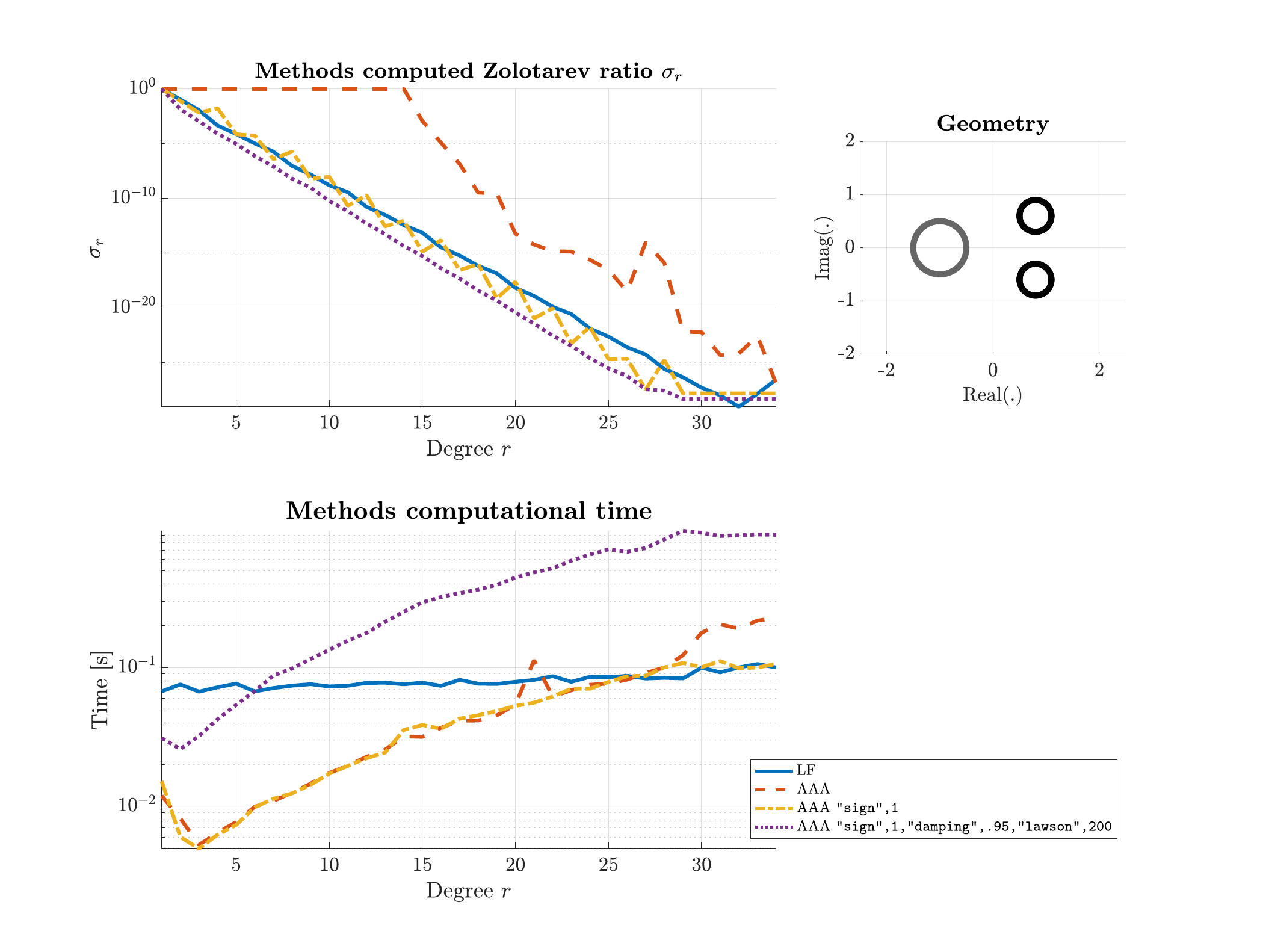}
    \caption{Case \texttt{'2a'}. Top: same description as \cref{fig:intro:1a}. Bottom: same description as \cref{fig:intro:1a_time}.}
    \label{fig:2a}
\end{figure}

This configuration exhibits symmetry with the  horizontal axis. Here again, the Zolotarev number obtained by the AAA sign and the AAA Lawson sign is most of the time lower than the LF. However, similarly to the previous observations and initial discussions, the LF is the only method able to reproduce a perfect complex conjugate collection of poles and zeros for both Z3 and Z4 approximations (see \cref{fig:2a}). This indicates that the structure of the true solution may be recovered by the LF. Note that the symmetry of poles and zeros means that the rational form can be obtained with real-valued polynomials only. Such a feature is relevant in the case where only data values are available, as it discovers the structure. 

\newpage
\subsection{Case \texttt{'2b'}}

\begin{figure}[H]
    \centering
    \includegraphics[width=.75\textwidth,align=c]{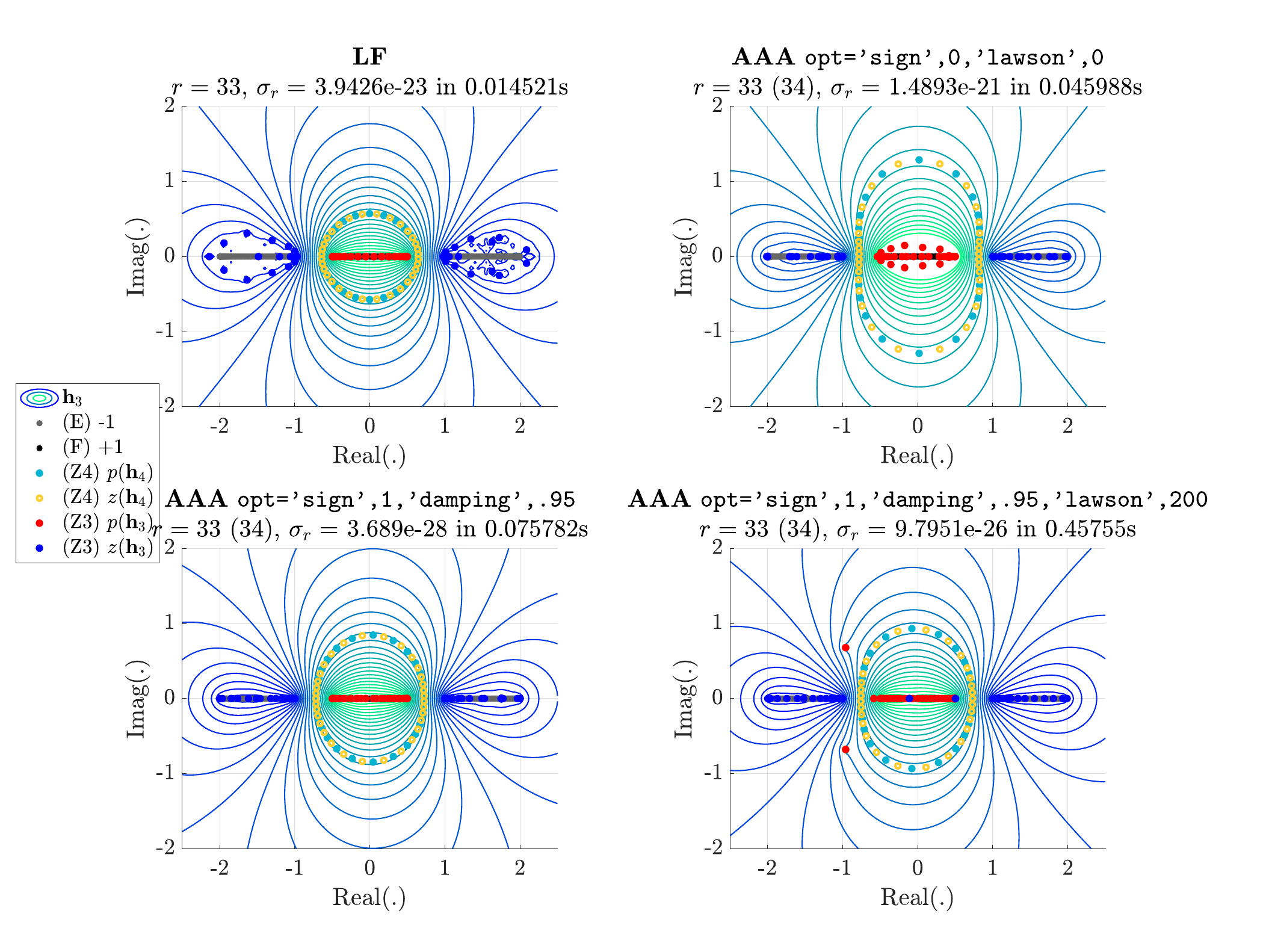}
    \includegraphics[width=.75\textwidth,align=c]{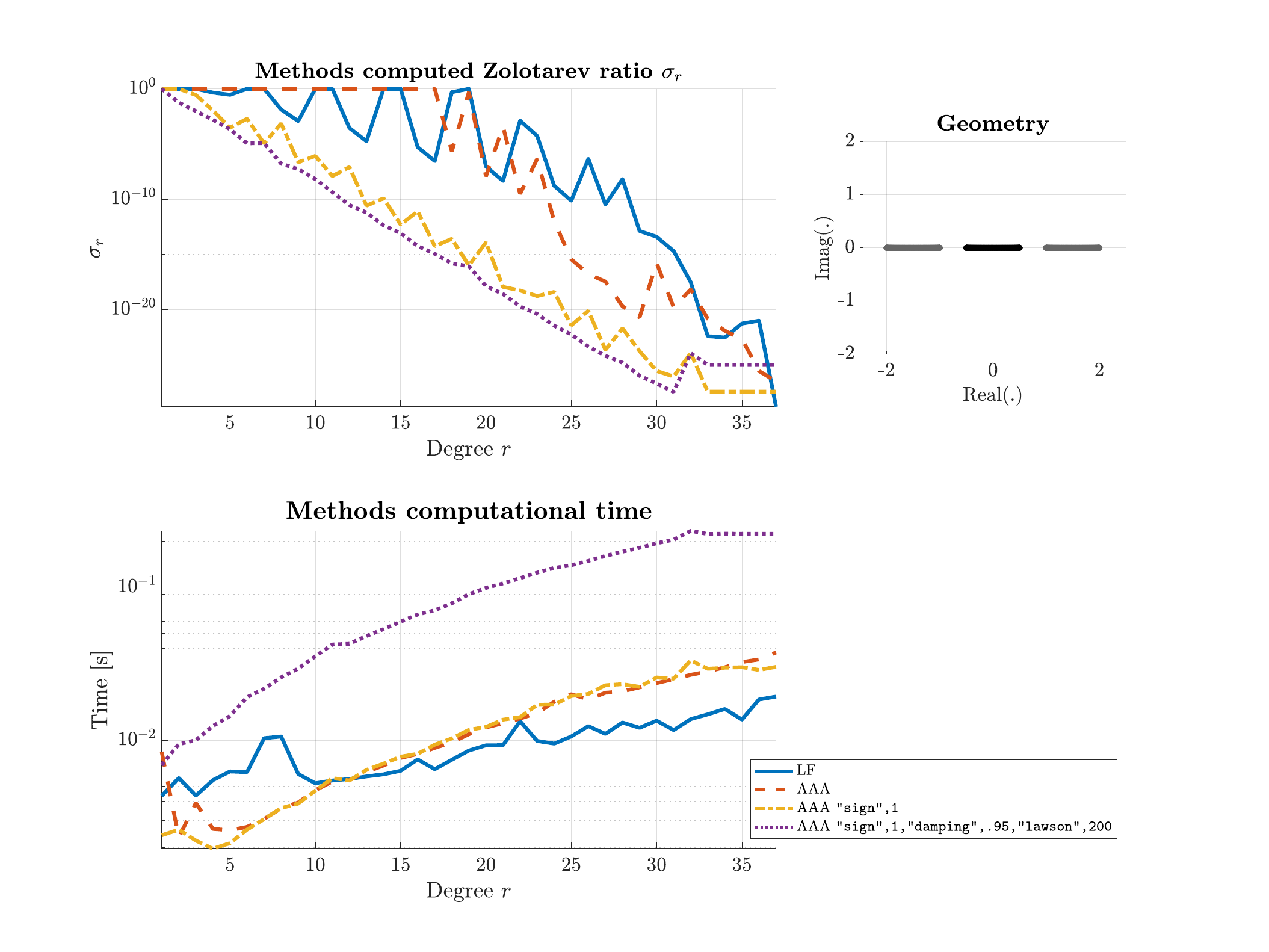}
    \caption{Case \texttt{'2b'}. Top: same description as \cref{fig:intro:1a}. Bottom: same description as \cref{fig:intro:1a_time}.}
    \label{fig:2b}
\end{figure}

This setup is an extension of the real-valued sign function \texttt{'1b'}, with two discontinuities. Same observations as the latter case can be made. The LF provides rather inaccurate results compared to the AAA sign and the AAA sign Lawson. However, the computation time is significantly lower. In this specific setup, the sign and Lawson variants clearly improve the accuracy. This observation is not surprising, as it gives substantial information to the approximation structure. However, by inspecting the upper frame of \cref{fig:2b}, the AAA sign Lawson approach results in poles and zeros with a strange  repartition, symptomatic of spurious poles/zeros distributions. In this upper frame, the LF does not recover a set of only real zeros for Z3; however, the latter are complex and closed by conjugation.

\newpage
\subsection{Case \texttt{'2c'}}

\begin{figure}[H]
    \centering
    \includegraphics[width=.75\textwidth,align=c]{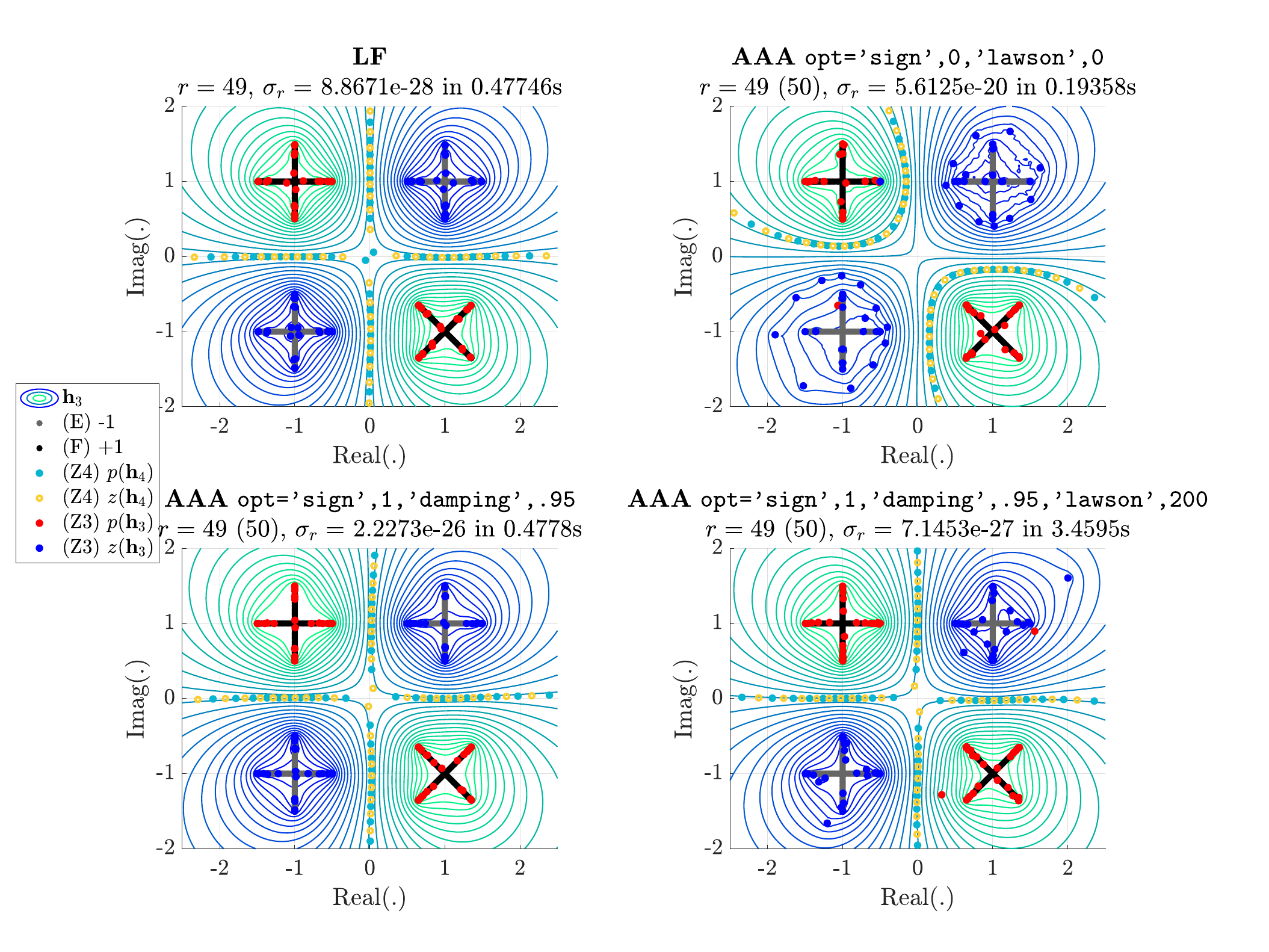}
    \includegraphics[width=.75\textwidth,align=c]{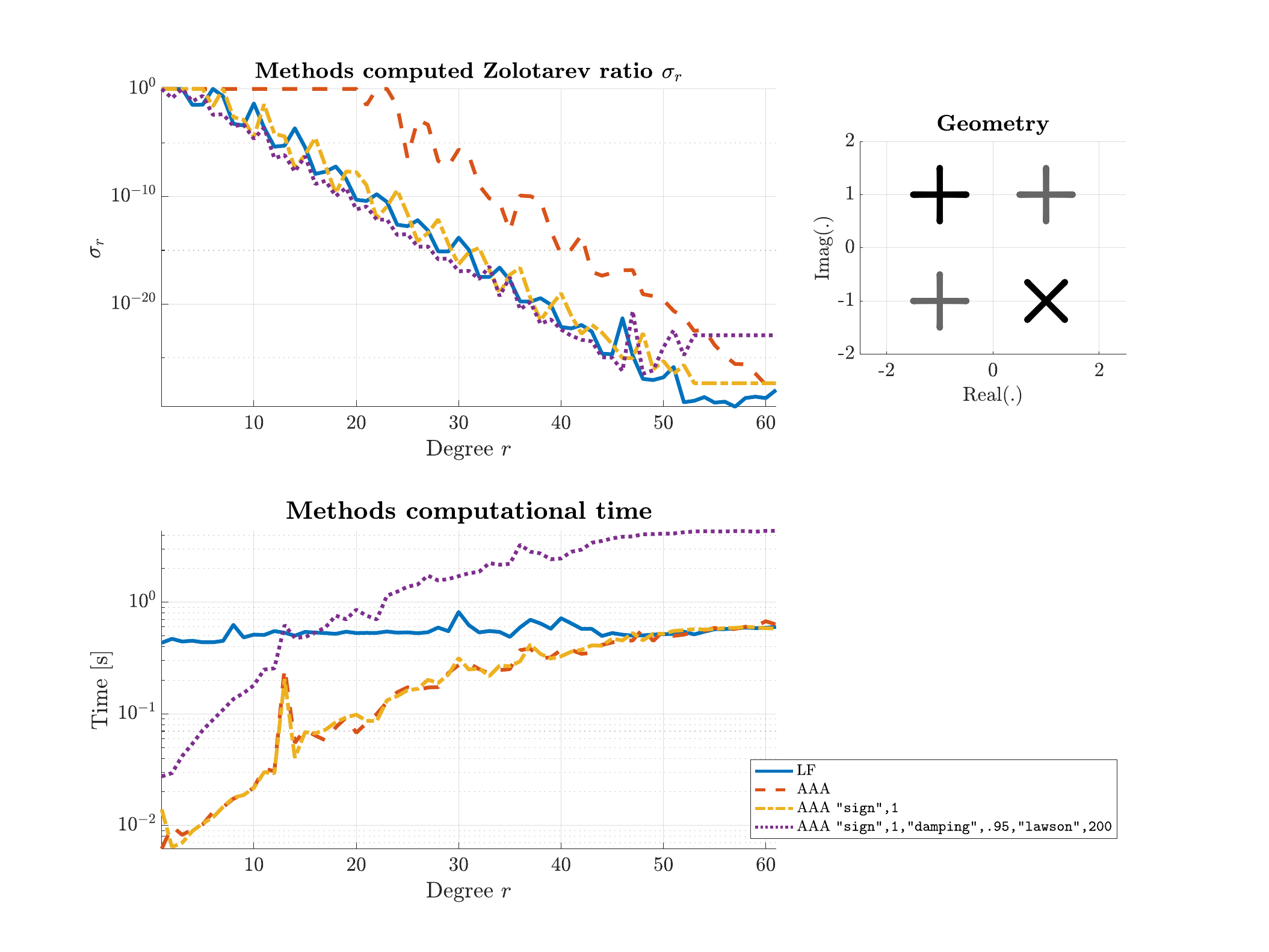}
    \caption{Case \texttt{'2c'}. Top: same description as \cref{fig:intro:1a}. Bottom: same description as \cref{fig:intro:1a_time}.}
    \label{fig:2c}
\end{figure}

This non-symmetric case illustrates a configuration where LF outperforms all other methods in both accuracy (see lower frame in \cref{fig:2c}). Once again, the AAA sign Lawson setup leads Z3 approximation  to poles and zeros peculiarly located with singularities and zeros in the non-expected area (top right and bottom right).

\newpage
\subsection{Case \texttt{'2d'}}

\begin{figure}[H]
    \centering
    \includegraphics[width=.75\textwidth,align=c]{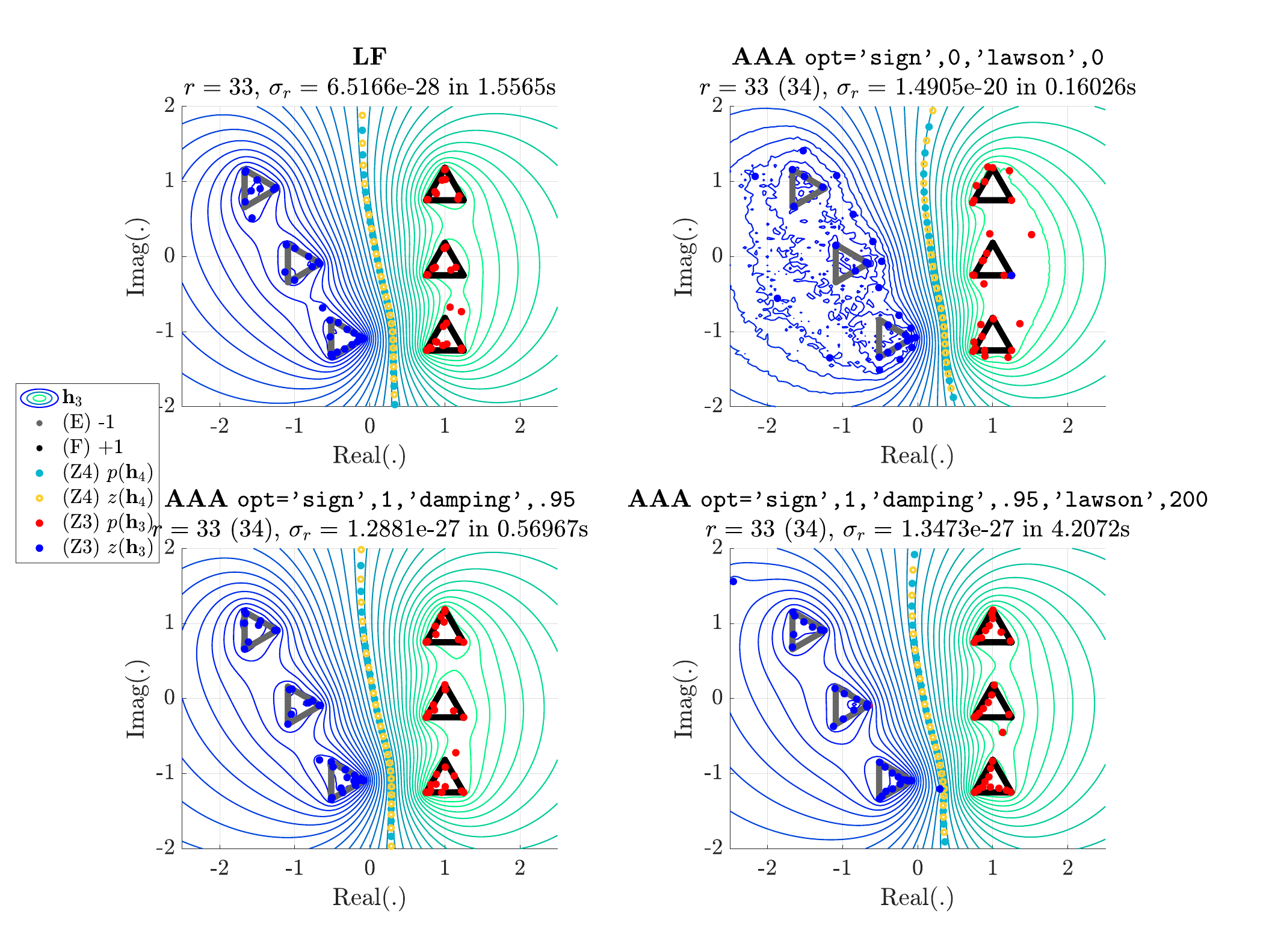}
    \includegraphics[width=.75\textwidth,align=c]{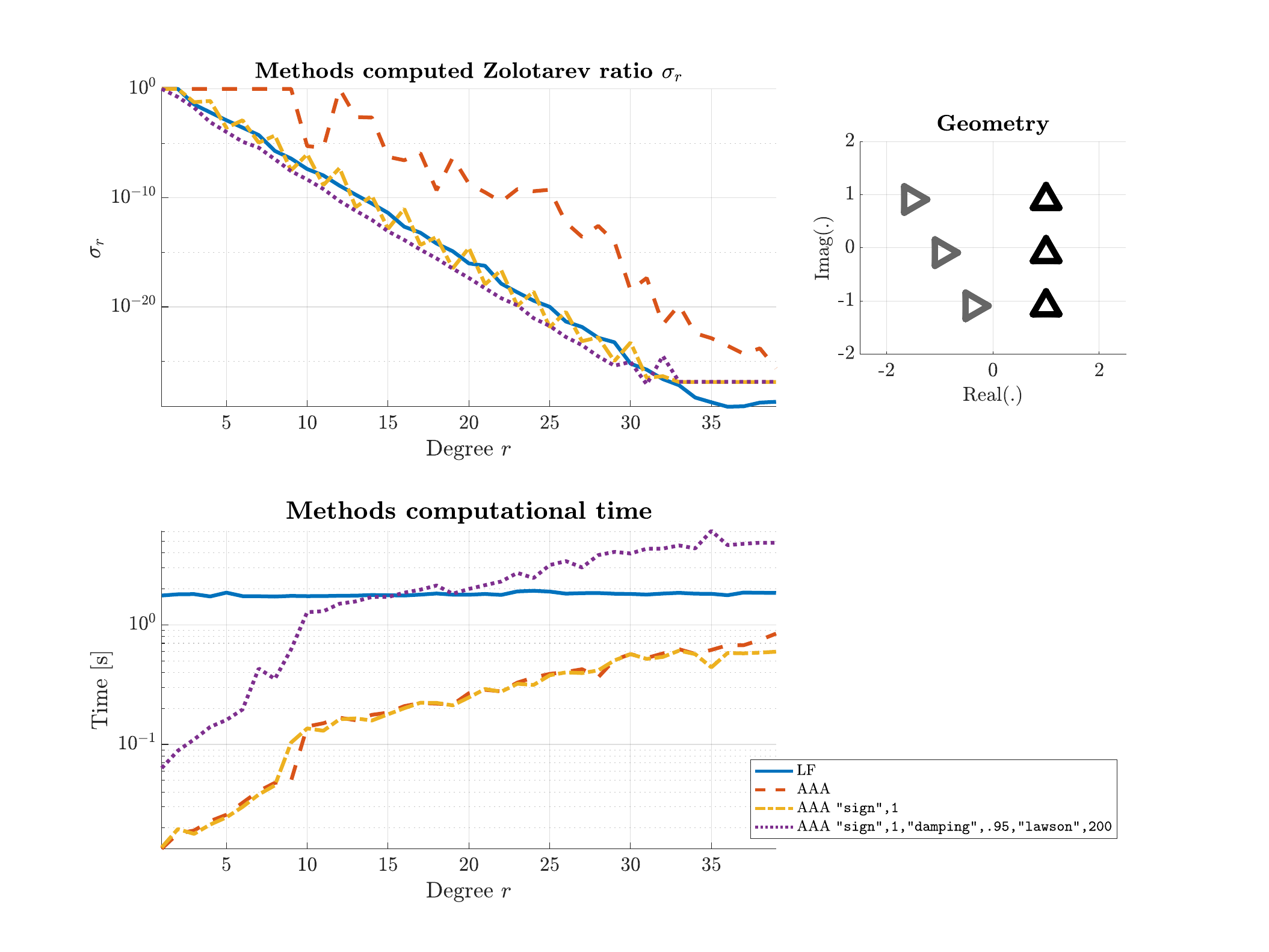}
    \caption{Case \texttt{'2d'}. Top: same description as \cref{fig:intro:1a}. Bottom: same description as \cref{fig:intro:1a_time}.}
    \label{fig:2d}
\end{figure}

This non-symmetric case illustrates a configuration where LF outperforms AAA, AAA with sign in accuracy (see lower frame in \cref{fig:2d}). In addition, the LF also provides better results than AAA with the sign and Lawson options when high accuracy is sought. 

\newpage
\subsection{Case \texttt{'3a'}, \texttt{'3b'}, \texttt{'3c'}, \texttt{'3d'}}

\begin{figure}[H]
    \centering
    \includegraphics[width=.45\textwidth,align=c]{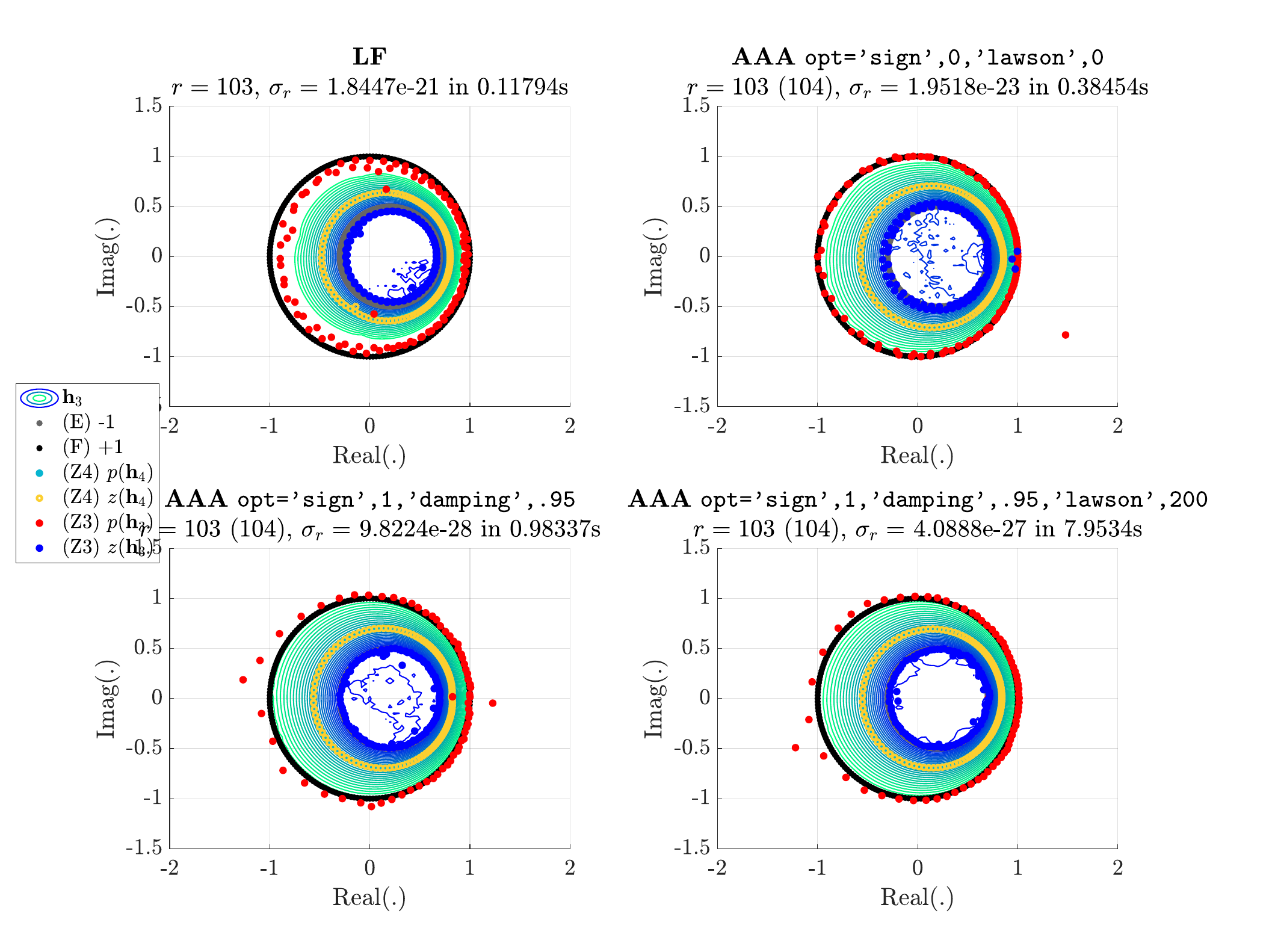}
    \includegraphics[width=.45\textwidth,align=c]{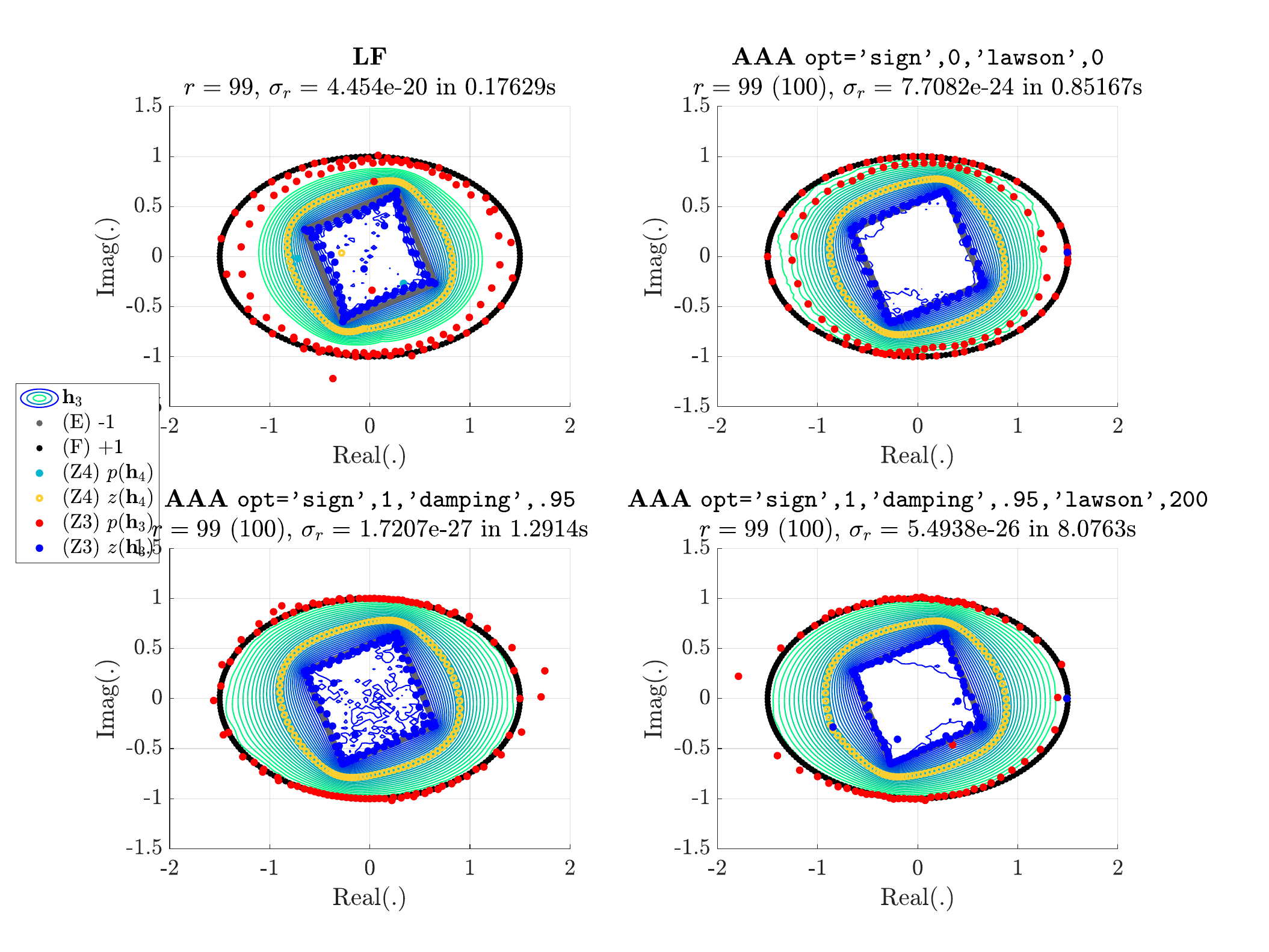}\\
    \includegraphics[width=.45\textwidth,align=c]{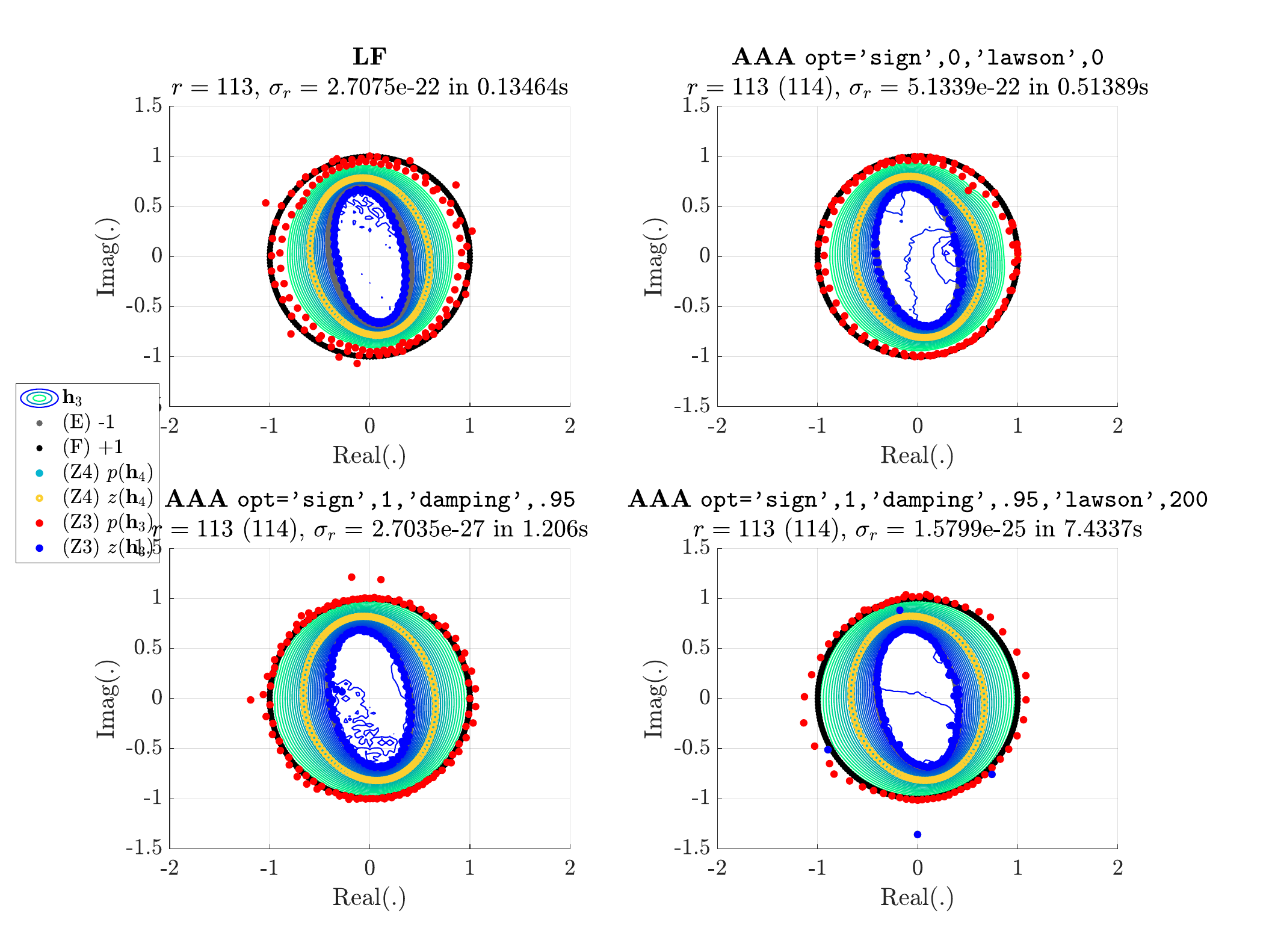}
    \includegraphics[width=.45\textwidth,align=c]{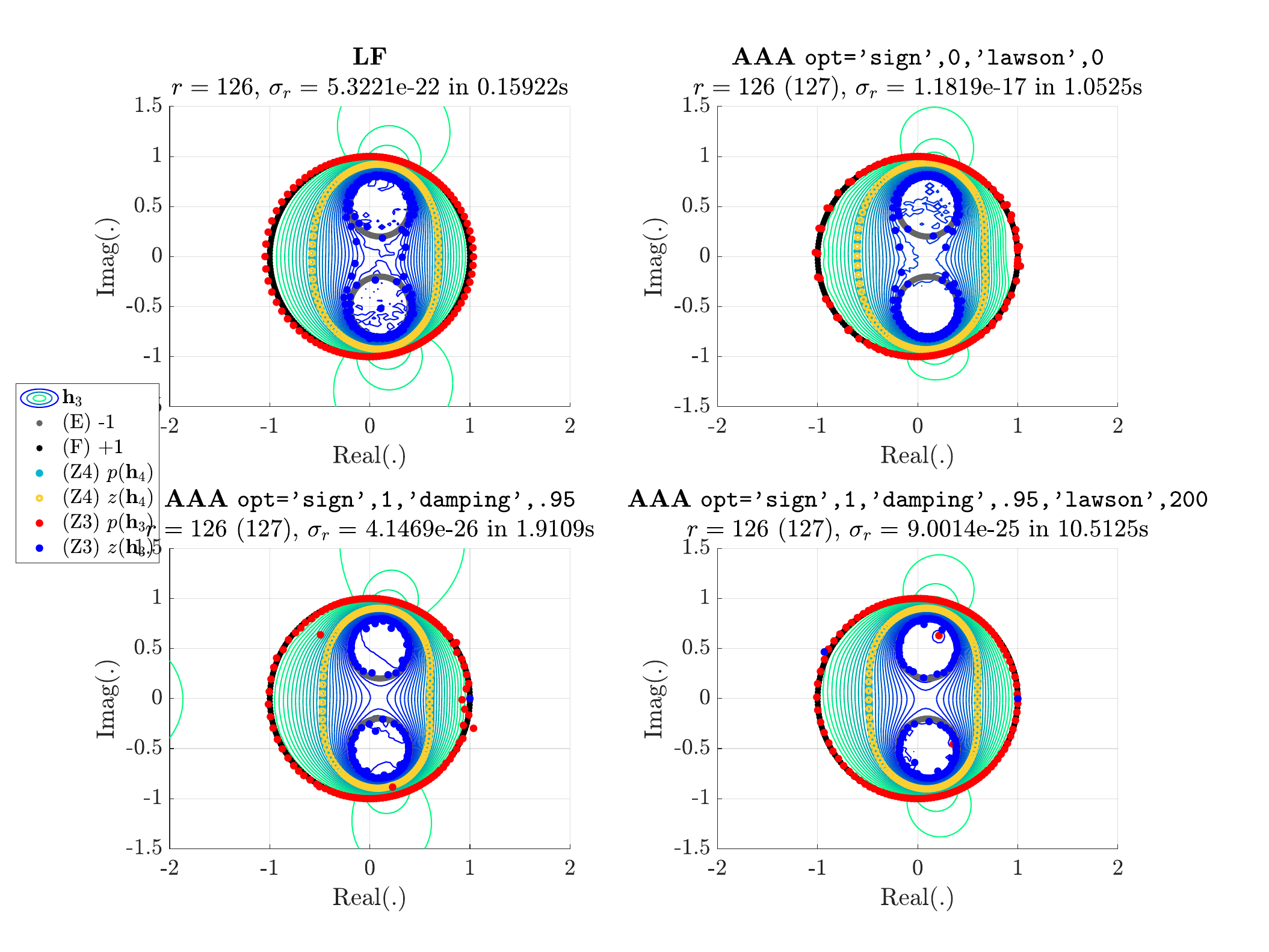}
    \caption{Case \texttt{'3a,b,c,d'}. Same description as \cref{fig:intro:1a}.}
    \label{fig:3a}
\end{figure}

\begin{figure}[H]
    \centering
    \includegraphics[width=.45\textwidth,align=c]{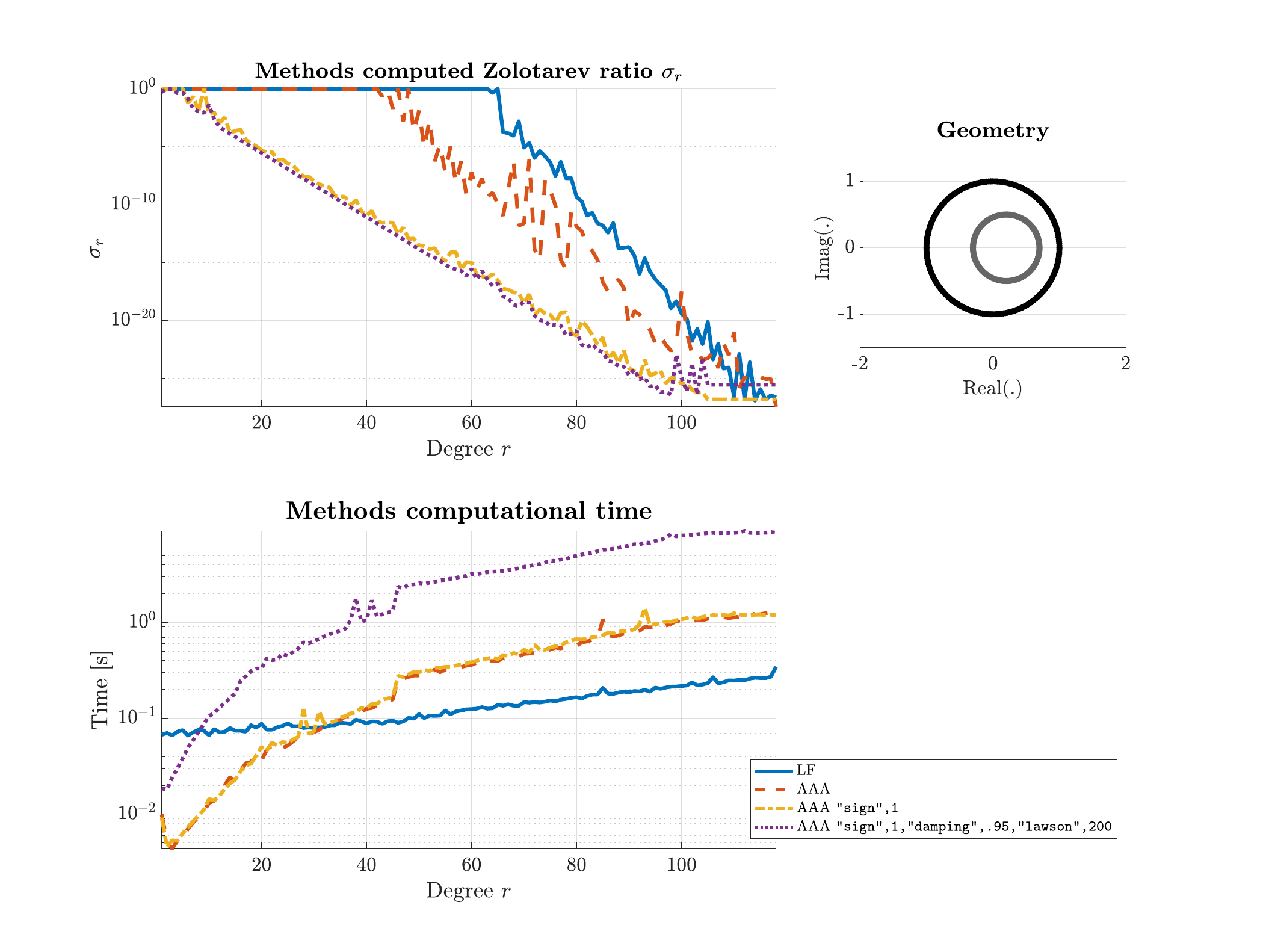}
    \includegraphics[width=.45\textwidth,align=c]{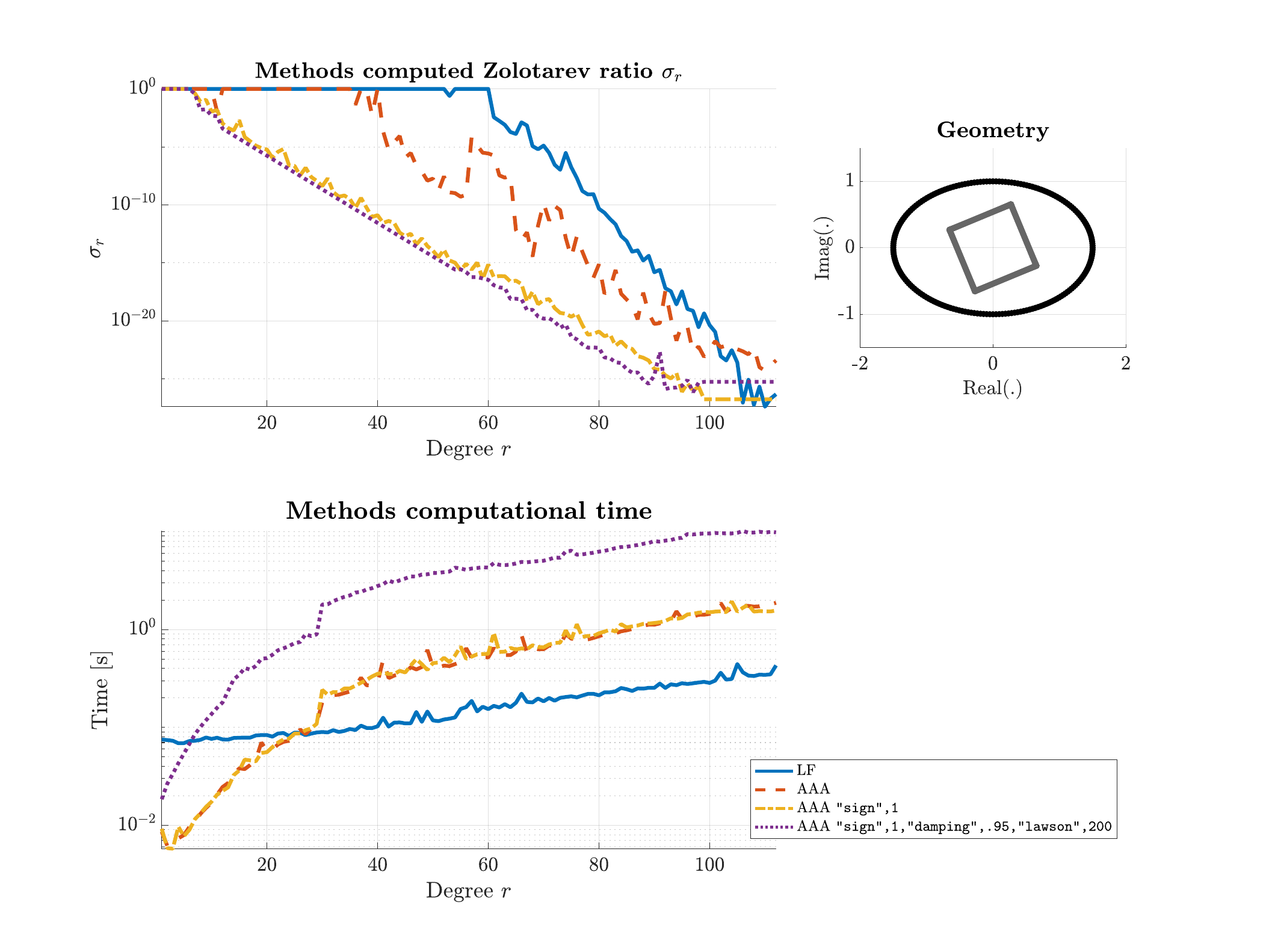}
    \includegraphics[width=.45\textwidth,align=c]{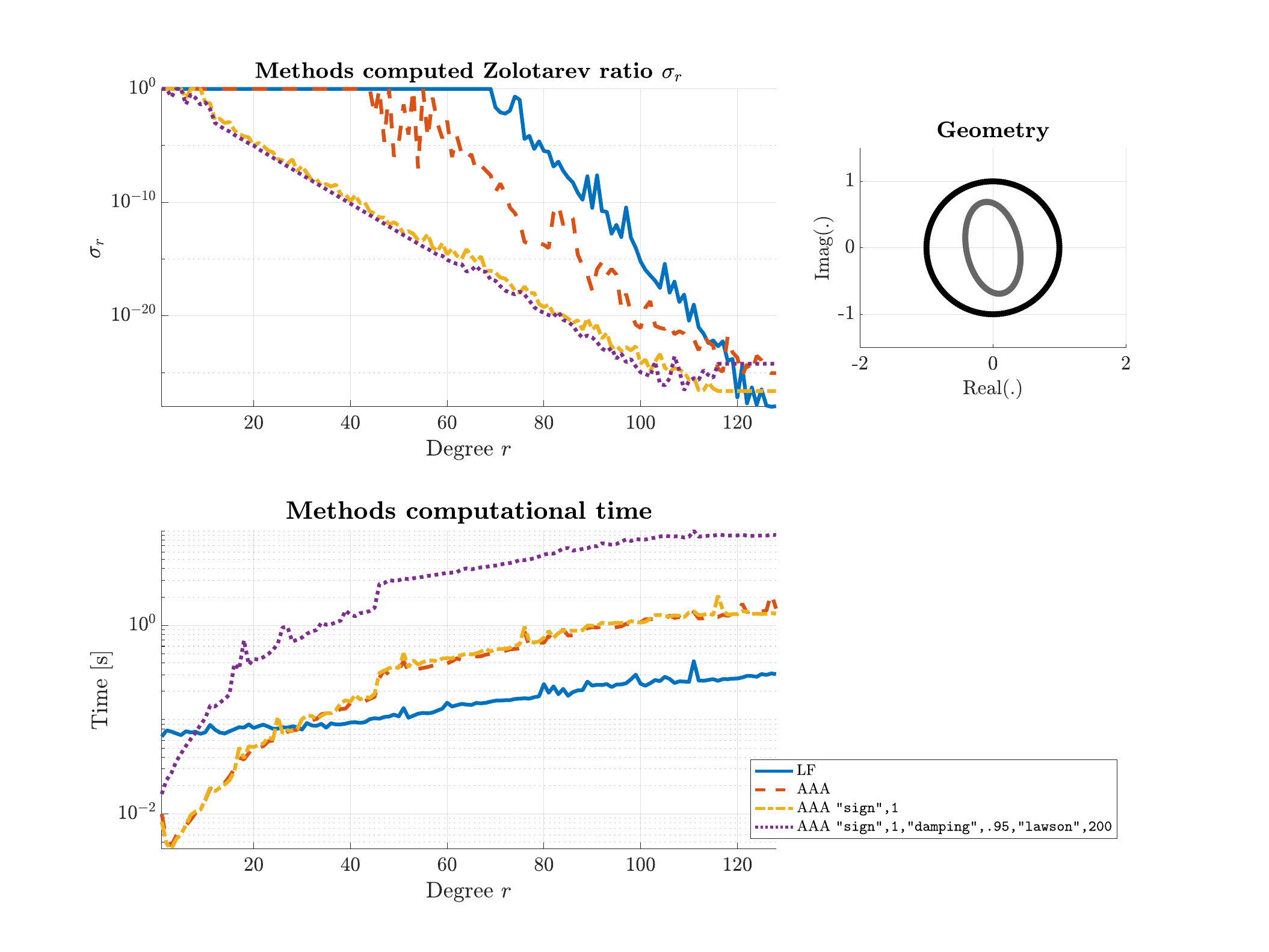}
    \includegraphics[width=.45\textwidth,align=c]{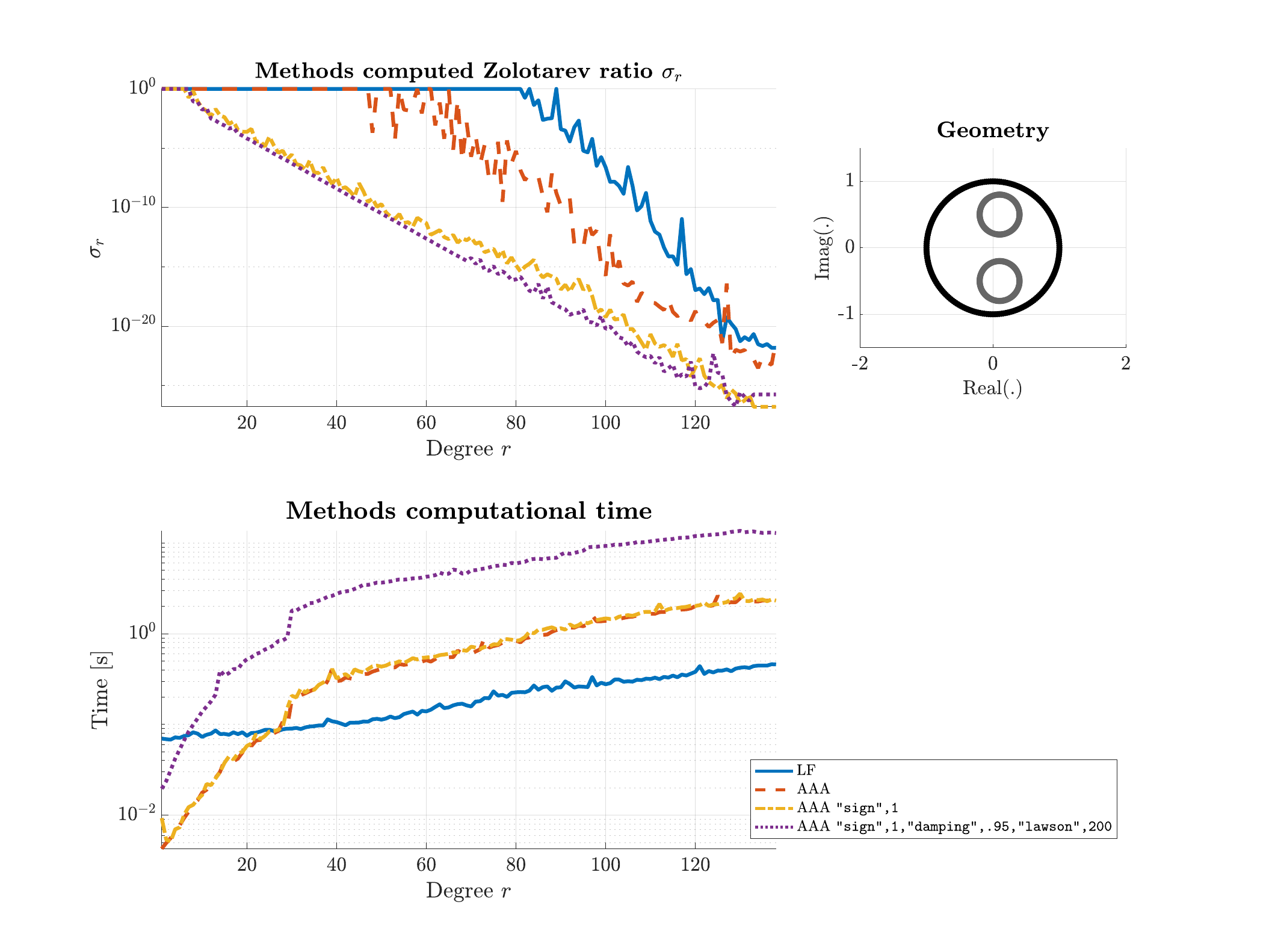}
    \caption{Case \texttt{'3a,b,c,d'}. Same description as \cref{fig:intro:1a_time}.}
    \label{fig:3b}
\end{figure}

Now we consider four configurations for which the sets are not disjoint, but $E$ is inside $F$. The results are reported in \cref{fig:3a,fig:3b}. In all these settings, the AAA configuration performs better than LF (except in computational time). It appears in all these cases that the LF needs a very high rational order  to achieve a good approximation, while the AAA sign and the AAA sign Lawson perform well already with low approximation orders. This aspect is not well understood yet, but certainly shows the usefulness of the Lawson and sign features in certain cases.

\newpage
\subsection{Case \texttt{'7'}}

\begin{figure}[H]
    \centering
    \includegraphics[width=.75\textwidth,align=c]{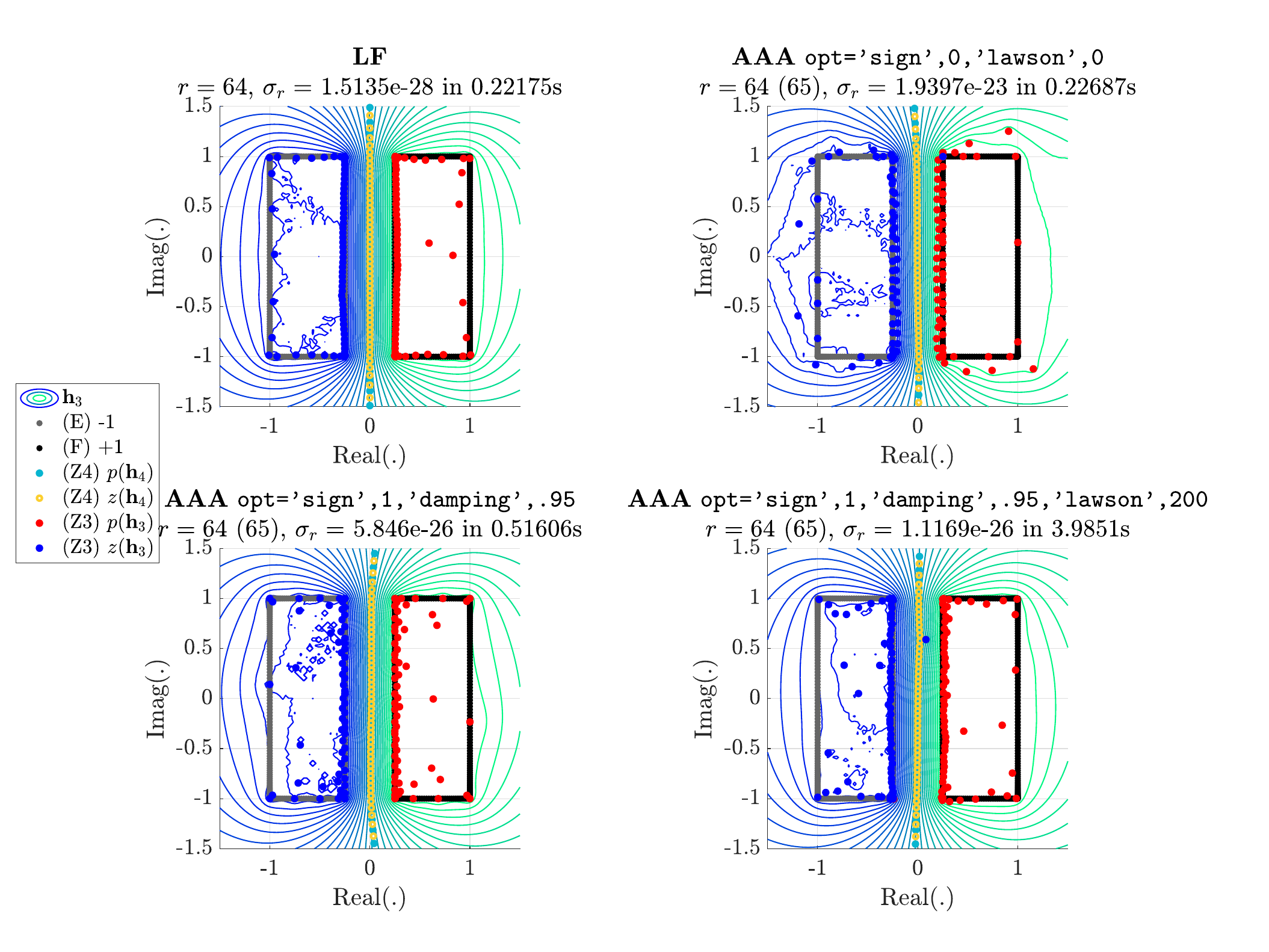}
    \includegraphics[width=.75\textwidth,align=c]{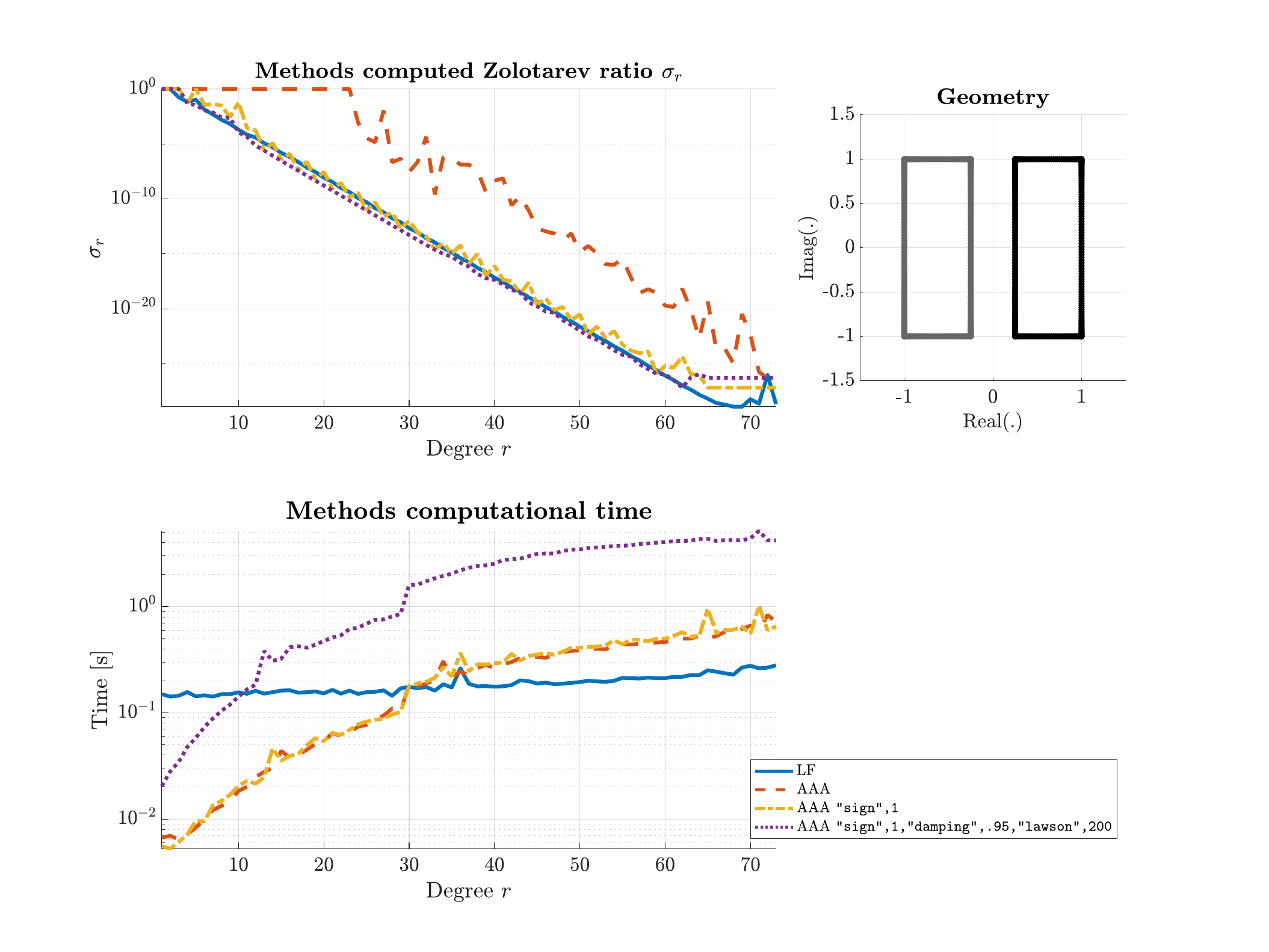}
    \caption{Case \texttt{'7'}. Top: same description as \cref{fig:intro:1a}. Bottom: same description as \cref{fig:intro:1a_time}.}
    \label{fig:7}
\end{figure}

The two-rectangle case shown in \cref{fig:7} illustrates the benefit of the LF with respect to AAA versions in terms of computational time. Again, the latter is perfectly flat and constant with almost the same accuracy. Moreover, in the configuration shown in the frame above, the approximation is even better with LF. Notice again the clean pole/zero distribution for the LF.

\newpage
\subsection{Case \texttt{'spiral1'}}

\begin{figure}[H]
    \centering
    \includegraphics[width=.75\textwidth,align=c]{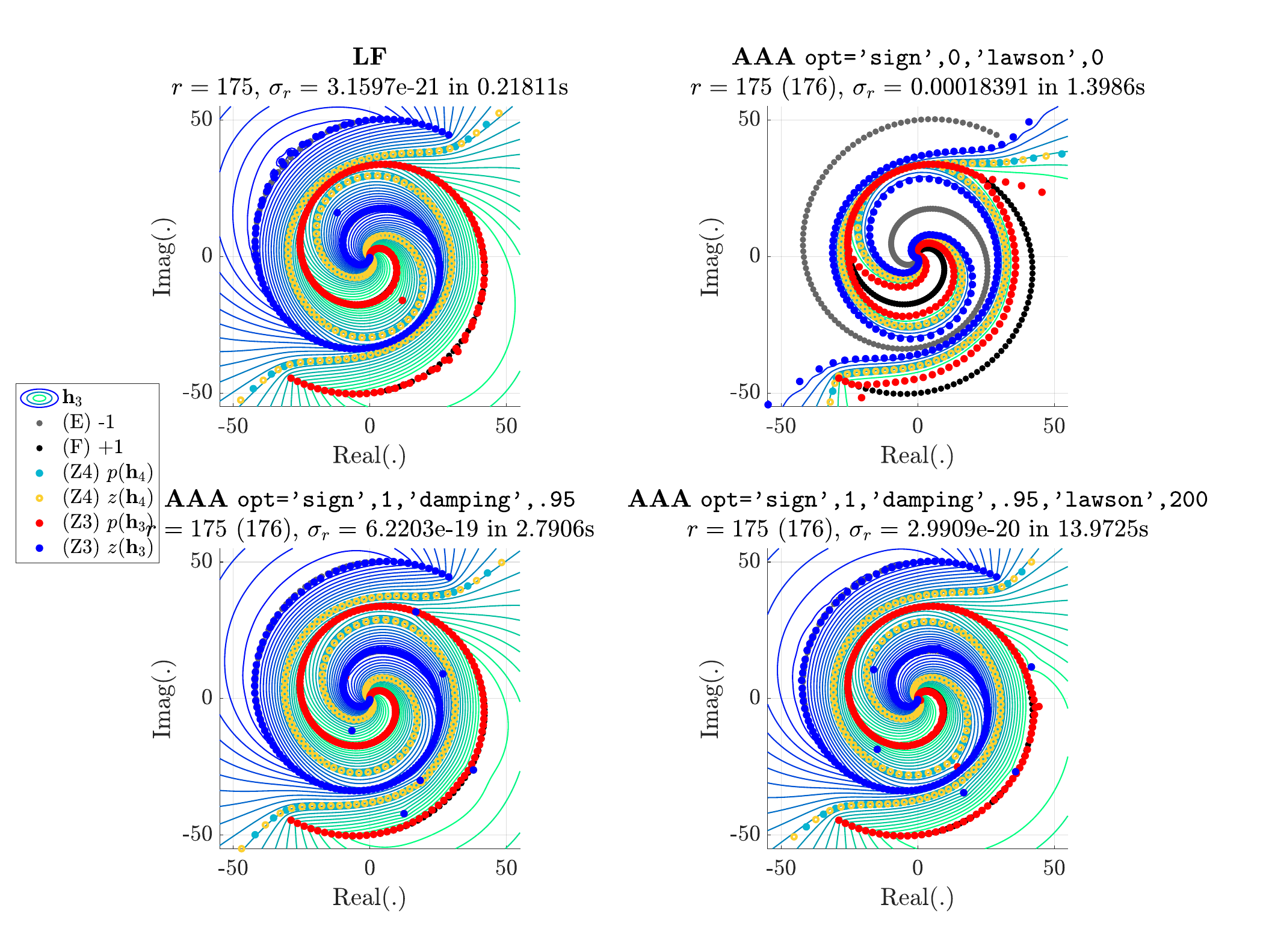}
    \includegraphics[width=.75\textwidth,align=c]{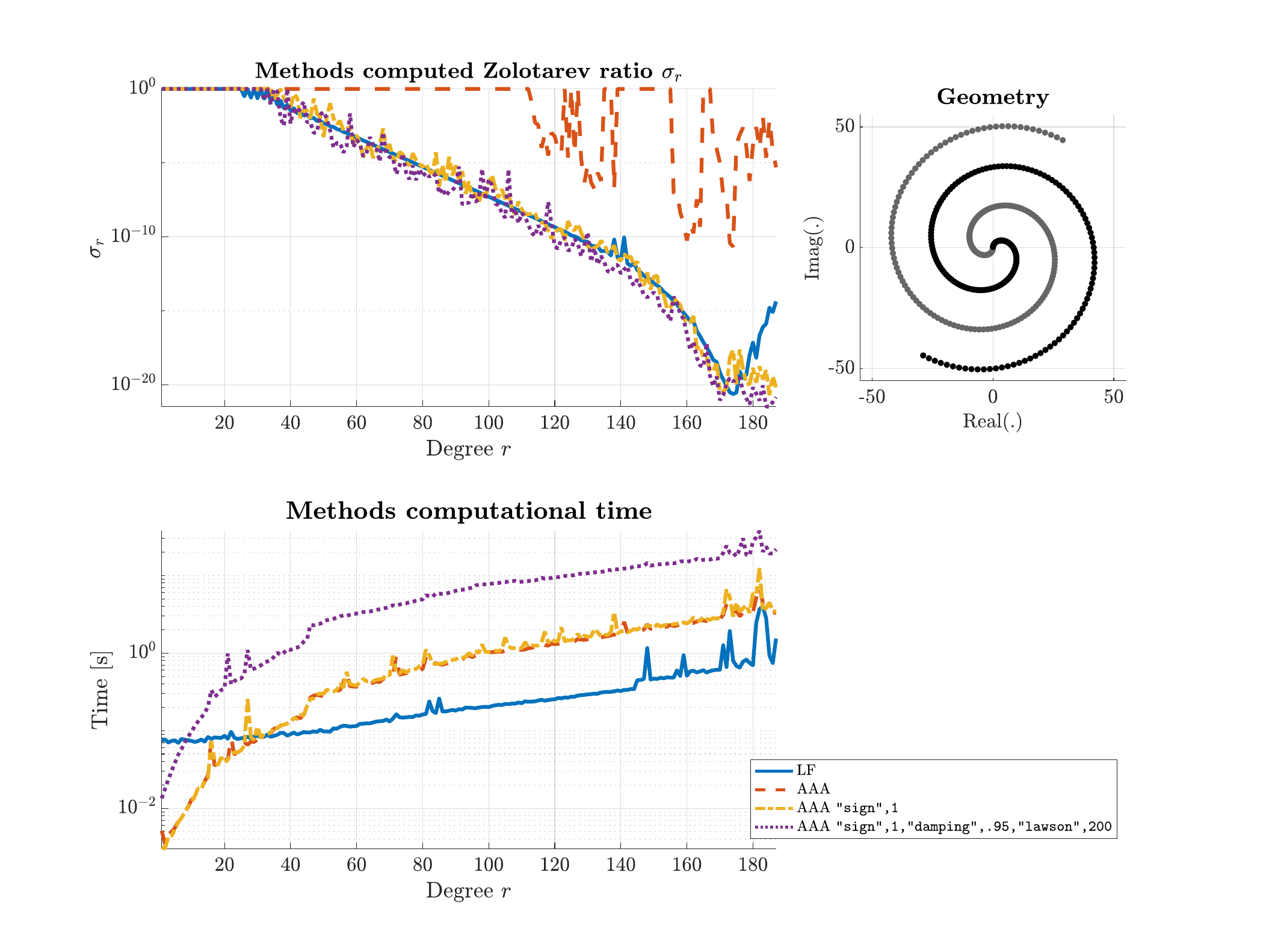}
    \caption{Case \texttt{'spiral1'}. Top: same description as \cref{fig:intro:1a}. Bottom: same description as \cref{fig:intro:1a_time}.}
    \label{fig:spiral1}
\end{figure}

The spiral example shown in \cref{fig:spiral1} (here computed with an LF normalized singular values threshold at $10^{-12}$) illustrates a configuration where again the LF outperforms in all situations, with a lower Zolotarev number, a lower computational time. Moreover, inspecting the upper part of \cref{fig:spiral1}, for the LF, only 2 poles and zeros seem out of the pattern, while the different AAA lead to a large number of irregularly spread poles and zeros. This lets us believe that the LF is able to recover the true structure of the optimal rational approximation.

\newpage
\subsection{Case \texttt{'pm2'}}

\begin{figure}[H]
    \centering
    \includegraphics[width=.75\textwidth,align=c]{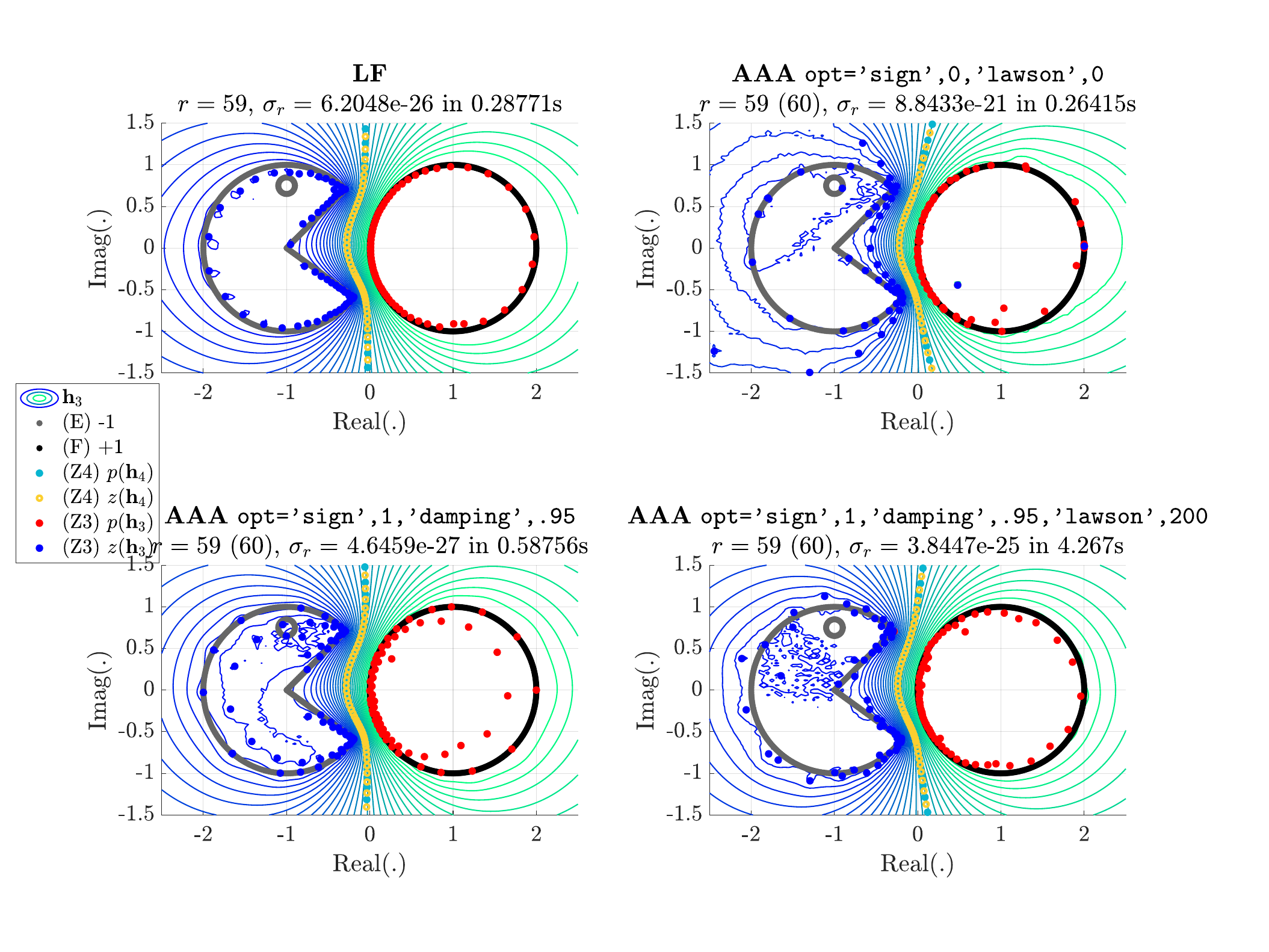}
    \includegraphics[width=.75\textwidth,align=c]{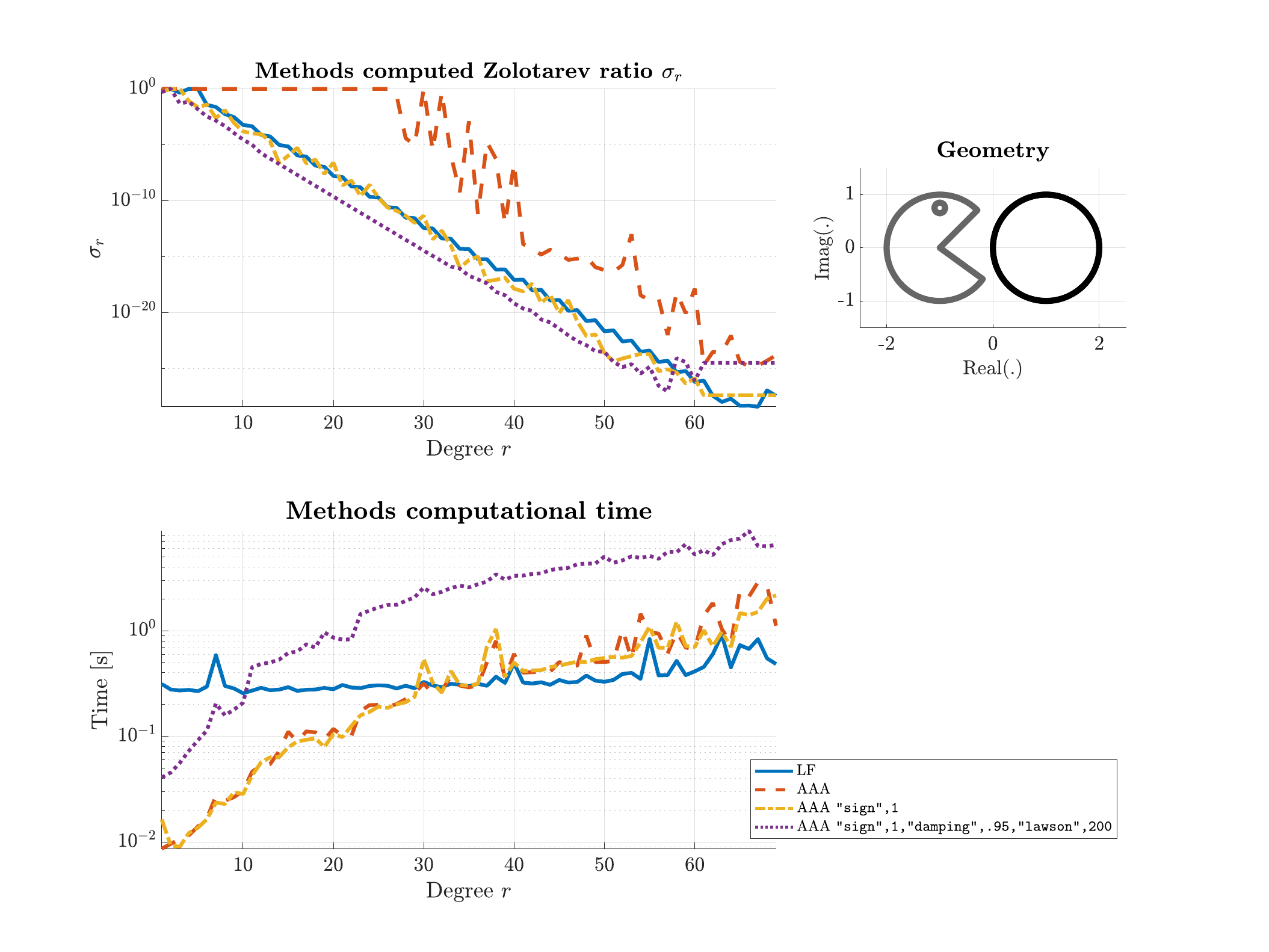}
    \caption{Case \texttt{'pm2'}. Top: same description as \cref{fig:intro:1a}. Bottom: same description as \cref{fig:intro:1a_time}.}
    \label{fig:pm2}
\end{figure}

The last configuration in \cref{fig:pm2} shows the Pac-Man figure (to the left) from the homonymous video game in the 80s \footnote{More details here: \url{https://pacman.com/en/history/}}, seeking to eat a candy depicted as a disk (on the right). Similar comments hold with almost the same accuracy between LF, AAA sign, and AAA sign Lawson, but with a significantly lower computation time for LF. Again, the nicely symmetric and well-distributed poles and zeros are obtained by the LF methodology only.

\section{Applications}\label{app:applications}
Application of the Zoloratarev problems are presented here. In  \cref{ssec:filter_iir}, we consider the infinite impulse response (IIR) filter design problem. 

\subsection{IIR digital filter design (Z4)}\label{ssec:filter_iir}

\begin{figure}[ht]
    \centering
    \includegraphics[width=1\textwidth]{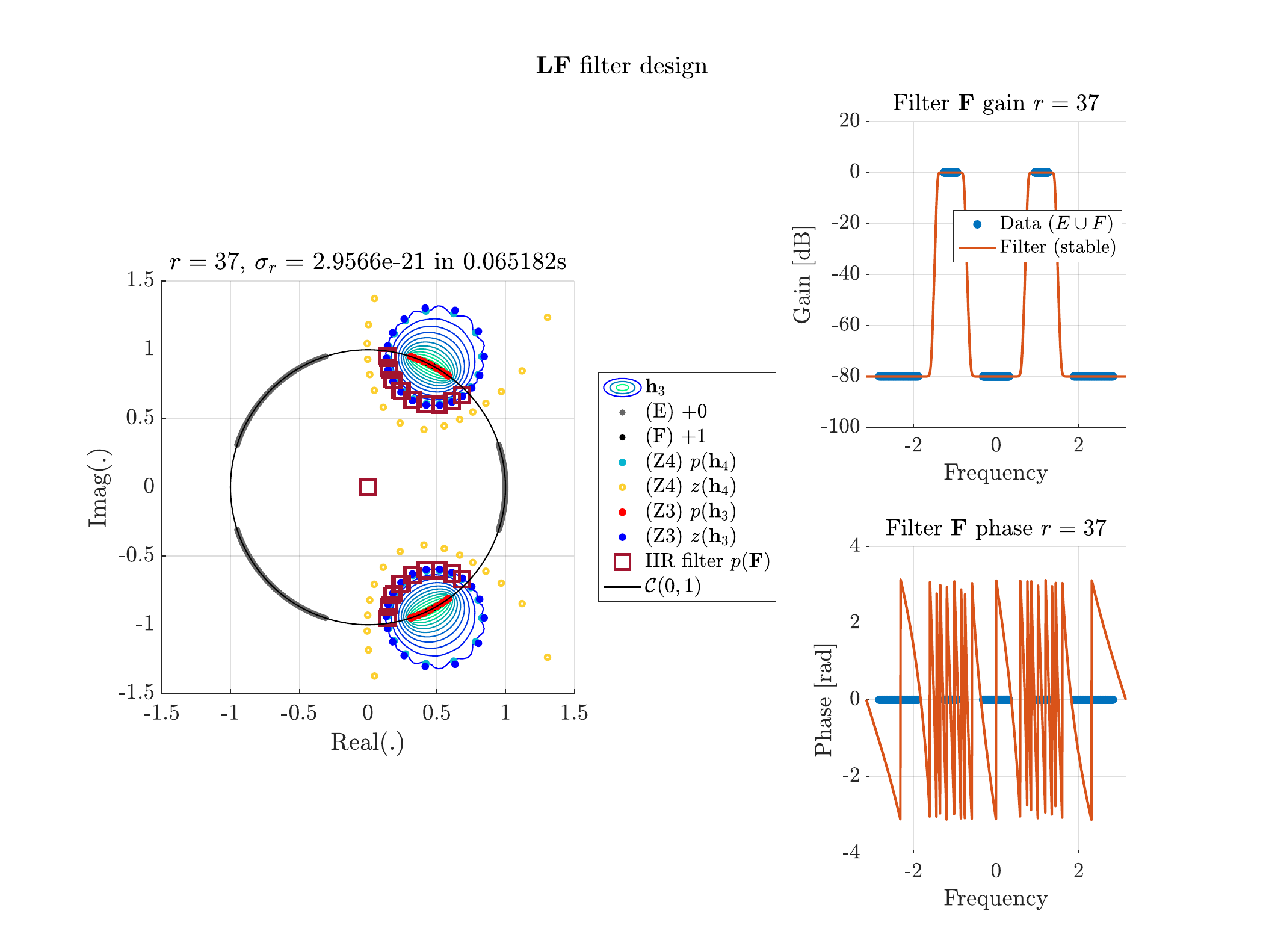}
    \caption{Bandpass filter and connection with Z4. Left: same description as \cref{fig:intro:1a} with additionally, the mirrored poles of Z4. Right: Bode gain and phase of the obtained filter $\mathbf{F}(\varIC)$.}
    \label{fig:filter}
\end{figure}

\begin{figure}[ht]
    \centering
    \includegraphics[width=1\textwidth]{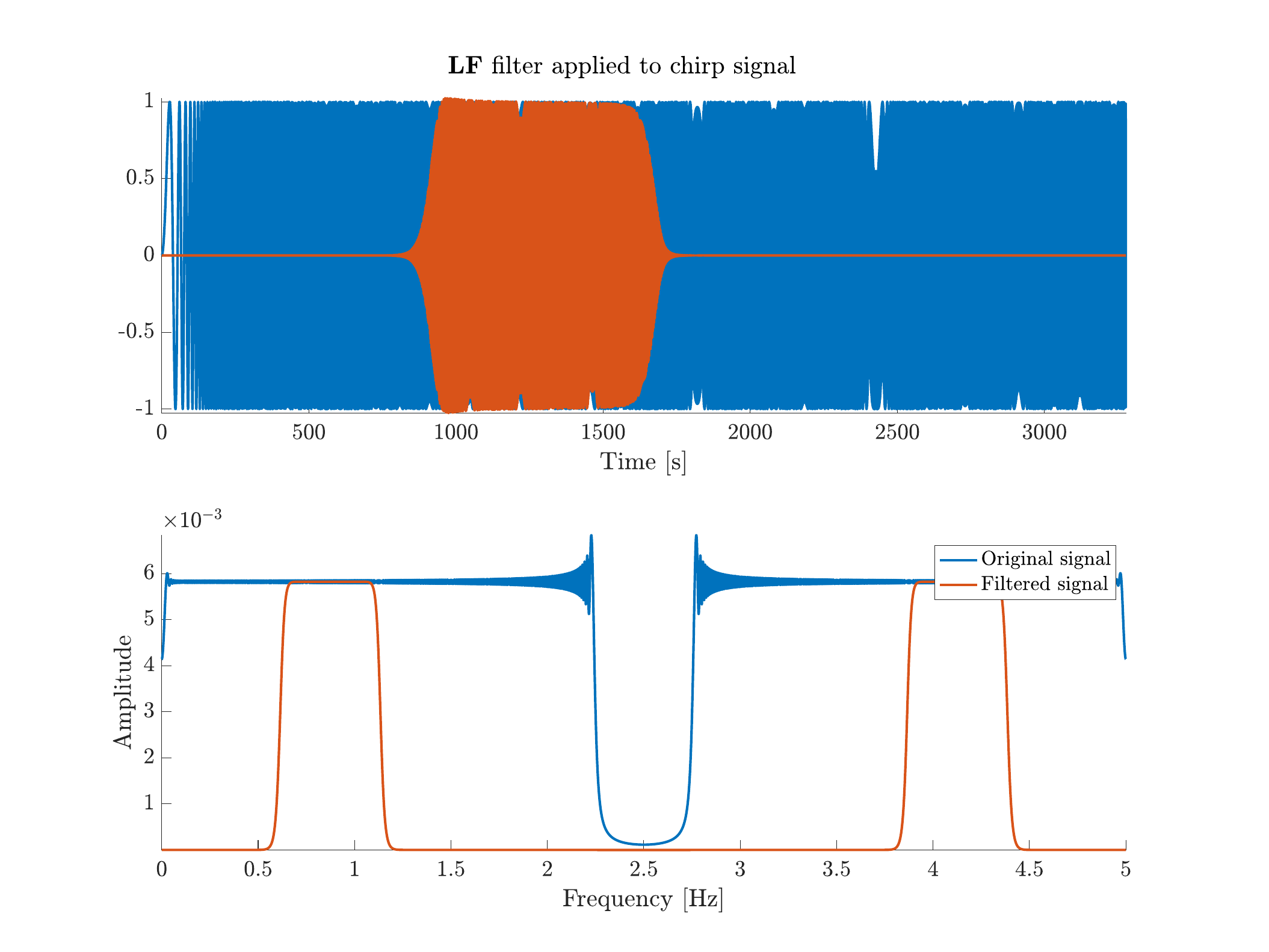}
    \caption{Application of filter $\mathbf F(\varIC)$ to a chirp signal. Top frame is the time-domain responses, bottom is the associated FFT absolute value.}
    \label{fig:filter_td}
\end{figure}

We are interested in designing a bandpass filter (for simplicity, we restrict the exposition to this class of filter but multi-band can easily be obtained). The problem consists in designing a filter which amplitude is $1$ over the frequencies of interest ($F$ set) and $\epsilon>0$ over the frequencies  to be discarded ($E$ set). This is precisely the Z4 sign problem. To apply the methodology presented, following steps are chosen\footnote{Filter and figures can be obtained with \texttt{demo3\_filterDesign.m}.}:
\begin{itemize}
\item Define $E$ ($=\epsilon$) and $F$ ($=1$) as
$$
\begin{array}{rcl}
E &=& e^{\imath \pi \Theta_1} \cup e^{\imath \pi \Theta_3} \\
F &=& e^{\imath \pi \Theta_2}
\end{array}
$$
with $\Theta_1=[10^{-3},0.1]$,  $\Theta_2=[0.3,0.4]$ and $\Theta_3=[0.6,0.9]$ sets. In the following application each set is linearly sampled between bounds with 50, 100 and 50 points. Refer to \cref{fig:filter} (left) to view the $E$ (grey dots) and $F$ (black dots) sets, and to \cref{fig:filter} (right) to view the corresponding values in Bode like form (blue dots).
\item  Then, apply the LF approach, as exposed in the main document. Here, as we seek for a digital filter that can be effectively implemented on a real-time system, real entries in the coefficients are expected only. Therefore, interpolation points $E\cup F$ are closed by conjugation before applying LF with the modification explained in \cite[\S 2.4.5]{ali17}. This is why on \cref{fig:filter} (right), data are both on positive and negative sides. The LF computes the order with a singular value threshold set to $10^{-12}$, and a very low  Z3 approximation metric $\sigma_r$ is reached, see  \cref{fig:filter} (left). Now the rational approximant $\mathbf h_4$ has real-valued coefficients only and both poles and zeros come closed by conjugation.
\item As revealed by the \cref{fig:filter} (left), some poles $p(\mathbf h_4)$ of the rational approximant are outside the unit circle $\mathcal C(0,1)$ (solid black line), which results in an unstable rational function, not implementable. Here, we chose to project the unstable poles inside the unit circle to preserve the gain. However, this indirectly impacts the phase in an uncontrollable manner, as shown in \cref{fig:filter} (right-bottom). 
\item As the coefficients are real, both the transfer function and realization are real. It is finally possible to create dynamical system with sampling time $T_s$ and to run it in real-time (see \cref{fig:filter_td}).
\end{itemize}

To convince of the applicability and interest of this methodology, \cref{fig:filter_td} illustrates the application of the obtained filter $\mathbf F(\varIC)$, with a sampling time $T_s=0.2$ seconds, on a chirp signal exciting frequencies from 0 to $0.9/(2T_s)$. In \cref{fig:filter_td}, the upper part represents the time-domain response, wile the lower is the frequency-domain response obtained by FFT.

\begin{remark}[Open questions]
The above process seems very simple, yet efficient. However, some legitimate questions arise: 
\begin{itemize}
\item The projection to enforce filter $\mathbf F$ stability modifies the phase. An appropriate criterion and method should be applied to control this step. As an example, (i) a linear phase with (almost) constant time group (e.g. Bessel like), or (ii) a minimal phase variation, can be of interest.
\item The gain data fit is very good on the chosen example. However, some configurations may result with spikes in the gain at the non interpolated values, which is unlikely for filtering. Avoiding these spikes in the non-interpolated data is important and methods should be thought, e.g. (i) data completion, (ii) constraints on the derivative of $\bF(\varIC)$, or (iii) a Lawson like method for LF may be some ideas.
\end{itemize}
\end{remark}

\end{document}